\documentclass[12pt,reqno]{amsart}

\begin{document}

\title[Bijective proofs for Schur function identities]{Bijective
        proofs for Schur function identities which imply
        Dodgson's condensation formula and Pl\"ucker relations}

\begin{abstract}

%%%@@@ input from ABSTRACT
We present a ``method'' for bijective proofs for
determinant identities, which is based on translating determinants
to Schur functions by the Jacobi--Trudi identity. We illustrate
this ``method'' by generalizing a bijective construction (which was
first used by Goulden) to a class of Schur function identities, from
which we shall obtain bijective proofs for Dodgson's condensation
formula, Pl\"ucker relations and a recent identity
of the second author.

%%%@@@ end input from ABSTRACT
\end{abstract}

\author{Markus Fulmek\\Michael Kleber}
\address{
Institut f\"ur Mathematik der Universit\"at Wien\\
Strudlhofgasse 4, A-1090 Wien, Austria.\newline\leavevmode\indent
Massachusetts Institute of Technology\\
77 Massachusetts Avenue, Cambridge, MA 02139, USA.
}
\email{
{\tt Markus.Fulmek@Univie.Ac.At}\\
{\tt Kleber@Math.Mit.Edu}
}

\date{\today}

\rm

\maketitle

%%%@@@ input from parts/macros
\def\qp#1#2{[{#1}^{#2}]}
\def\sqp#1#2{s_{\qp{#1}{#2}}}
\def\pa#1{(#1)}
\def\sla#1{s_{\pa{#1}}}
\def\TP{P}
\def\TQ{A}
\def\TS{B}
\def\minor#1#2#3#4#5{#1_{\{#2,#3\},\{#4,#5\}}}
\def\la{\lambda}
\def\si{\sigma}
\def\pprime{{\prime\prime}}

%%%%%%%%%%%%%%%%%%%%%%%%%%%%%%%%%%%%%%%%%%%%%%%%%%%%%%%%%%%%%%%%%%%%%%%%%%%%
%                                       Diverse Kleinigkeiten, Abkuerzungen etc.:                                                       %
%%%%%%%%%%%%%%%%%%%%%%%%%%%%%%%%%%%%%%%%%%%%%%%%%%%%%%%%%%%%%%%%%%%%%%%%%%%%
\let\epsilon\varepsilon
\newcommand{\esf}{{e}}                  % elementary symmetric function
\newcommand{\hsf}{{h}}                  % homogeneous symmetric function
\newcommand{\ksf}{{k}}                  % Delta^2 of homogeneous symmetric function
\newcommand{\schf}[2]{{s}_{#1}\of{#2,\bold{x}}}
                                                                                                % Schur function
\newcommand{\defeq}{:=}                                 % define equal...
\newcommand{\tableau}{T}                                % semistandard tableau
\newcommand{\Ferrer}{B}                            % Ferrers board
\newcommand{\digraph}{D}                                % directed graph
\newcommand{\sdigraph}{S}                               % Stembridge's directed graph
\newcommand{\horedge}{a}                                % horizontal edge
\newcommand{\veredge}{\hat{a}}          % vertical edge
\newcommand{\path}{P}                                   % lattice path
\newcommand{\pathset}{M}                                % lattice path's adjoined edge set
\newcommand{\fpath}{{\mathcal P}}                       % family of lattice paths
\newcommand{\nfpath}{{\cal N}}          % number of lattice paths
\newcommand{\nifpath}{{\cal P}_o}       % nonintersecting family of lattice paths
\newcommand{\SG}[1]{\boldsymbol{S}_{#1}}
                                                                                                % Symmetric group
\newcommand{\GF}[1]{\bold{GF}\of{#1}}
                                                                                                % Generating function
\newcommand{\sgn}{\operatorname{sgn}}                   % Signum
\newcommand{\inv}{\operatorname{inv}}                   % Inversionen
\newcommand{\sschf}[2]{{sp}_{#1}\of{#2,\bold{x}}}
\newcommand{\isschf}[2]{{sp}_{n,m}\of{#1,#2;\bold{x},\bold z}}
\newcommand{\osschf}[2]{{sp}_{n,1}\of{#1,#2,\bold{x},z}}
                                                                                                % symplectic Schur function
\newcommand{\reflect}{{\bold{R}}}               % reflection
\newcommand{\mreflect}{\tilde{\bold{R}}}                % modified reflection
\newcommand{\oschf}[2]{{o}_{#1}\of{#2,\bold{x}}}
                                                                                                % orthogonal Schur function
\newcommand{\oschftriv}[2]{{o}_{#1}\of{#2,\seqof{1,1,\dots,1}}}
                                                                                                % "coarse" orthogonal Schur function
\newcommand{\liff}{\text{ iff }}                % logical: if and only if
\newcommand{\limp}{\text{ implies }}
                                                                                                % logical: implies
\renewcommand{\land}{\text{ and }}      % logical: and

\newcommand{\thmref}[1]{Theorem~\ref{#1}}
\newcommand{\secref}[1]{\S\ref{#1}}
\newcommand{\lemref}[1]{Lemma~\ref{#1}}
\newcommand{\figref}[1]{Figure~\ref{#1}}

%%%%%%%%%%%%%%%%%%%%%%%%%%%%%%%%%%%%%%%%%%%%%%%%%%%%%%%%%%%%%%%%%%%%%%%%%%%%
%                    Diverse Klammerungen:                                 %
%%%%%%%%%%%%%%%%%%%%%%%%%%%%%%%%%%%%%%%%%%%%%%%%%%%%%%%%%%%%%%%%%%%%%%%%%%%%
\newcommand{\seqof}[1]{\left\langle#1\right\rangle}
\newcommand{\setof}[1]{\left\{#1\right\}}
\newcommand{\of}[1]{\left(#1\right)}
\newcommand{\parof}[1]{\left(#1\right)}
\newcommand{\numof}[1]{\left|#1\right|}
\newcommand{\detof}[2]{\left|{#2}\right|_{#1\times #1}}
\newcommand{\firstcolumn}[1]{%\underbrace{#1%}_{\text{first column}}
\quad \vdots\quad }
\newcommand{\firstrow}[1]{%\underbrace{#1%}_{\text{first row}}
\quad \vdots\quad }
\newcommand{\brkof}[1]{\left[#1\right]}
\newcommand{\intof}[1]{\lfloor #1\rfloor}
\newcommand{\absof}[1]{\left|#1\right|}
\newcommand{\sgnof}[1]{\operatorname{sgn}\of{#1}}

\newtheorem{thm}{Theorem}
\newtheorem{lem}[thm]{Lemma}
\newtheorem{dfn}[thm]{Definition}
\newtheorem{obs}[thm]{Observation}
\newtheorem{rem}[thm]{Remark}

%%%@@@ end input from parts/macros

%%%@@@ input from paths
%%%%%%%%%%%%%%%%%%%%%%%%%%%%%%%%%%%%%%%%%%%%%%%%%%%%%%%%%%%%%%%%%%%%%%%%%%%%
%               PATHS.TEX: Begonnen am 16. Juli 1993                       %
% Dieses kleine Macropaket soll einfache Befehle zum Zeichnen von Gittern  %
% und Gitterpunktwegen einfuehren. (Fuer AMSLaTeX!)                        %
% Wesentlich ueberarbeitet und gestrafft am 1. Mai 1994.                   %
% Ergaenzt am 24. August 1994: Autmoatisches labelling von Wegen,          %
% Giambelli-Wege ...                                                       %
%%%%%%%%%%%%%%%%%%%%%%%%%%%%%%%%%%%%%%%%%%%%%%%%%%%%%%%%%%%%%%%%%%%%%%%%%%%%

%%%%%%%%%%%%%%%%%%%%%%%%%%%%%%%%%%%%%%%%%%%%%%%%%%%%%%%%%%%%%%%%%%%%%%%%%%%%
% Fonts:                                                                   %
%%%%%%%%%%%%%%%%%%%%%%%%%%%%%%%%%%%%%%%%%%%%%%%%%%%%%%%%%%%%%%%%%%%%%%%%%%%%
\font\scalefont = cmti8
 
%%%%%%%%%%%%%%%%%%%%%%%%%%%%%%%%%%%%%%%%%%%%%%%%%%%%%%%%%%%%%%%%%%%%%%%%%%%%
% Basic parameters:                                                        %
%%%%%%%%%%%%%%%%%%%%%%%%%%%%%%%%%%%%%%%%%%%%%%%%%%%%%%%%%%%%%%%%%%%%%%%%%%%%
\unitlength = 5mm     % instead of 4 mm
\thicklines

\newcount\xmin                % coordinates of extremal points of the grid
\newcount\xmax
\newcount\ymin
\newcount\ymax

\newcount\gridwidth
\newcount\gridheight

\newcount\nofxpoints          % number of points in horizontal direction
\newcount\nofypoints          % number of points in vertical direction

\newcount\dcnta               % dummy-counter A
\newcount\dcntb               % dummy-counter B
\newcount\dcntc     % dummy-counter C

% Special macros for drawing a grid with lower left corner (xmin, ymin)
% and upper right corner (xmax, ymax):
\def\begingrid#1#2#3#4{   % arguments: xmin, ymin, xmax, ymax
 \global\xmin = #1
 \global\ymin = #2
 \global\xmax = #3
 \global\ymax = #4
 % simple test:
 \ifnum\xmin > \xmax\errmessage{PATHS: \xmin > \xmax|}\fi
 \ifnum\ymin > \ymax\errmessage{PATHS: \ymin > \ymax|}\fi
 % computation of heigth and width
 \global\gridwidth = \xmax
 \global\gridheight = \ymax
 \global\advance\gridwidth by -\xmin
 \global\advance\gridheight by -\ymin
 \nofxpoints = \gridwidth
 \advance\nofxpoints by 1
 \nofypoints = \gridheight
 \advance\nofypoints by 1
  % adjust some space around the picture
  % \advance\gridwidth by 4
  % \advance\gridheight by 4
  % \advance\xmin by -2
  % \advance\ymin by -2
 % LaTeX's picture command
 \begin{picture}(\gridwidth, \gridheight)(\xmin, \ymin)
  % restore correct the lengths
  % \advance\gridwidth by -4
  % \advance\gridheight by -4
  % \advance\xmin by 2
  % \advance\ymin by 2
 % draw gridpoints:
 \global\dcnta = \ymin
 \loop\ifnum\dcnta<\ymax
  \makegridline\dcnta % do inner loop in a separate macro
  \global\advance\dcnta by 1
 \repeat
 \makegridline\ymax
}
 
% inner loop of grid-drawing procedure
\def\makegridline#1{
   \begingroup
   \global\dcntb = \xmin
   \loop\ifnum\dcntb<\xmax
      \put(\dcntb, #1){\circle*{0.1}}
      \global\advance\dcntb by 1
   \repeat
   \put(\xmax, #1){\circle*{0.1}}
   \endgroup
}
 
% draw axes
\def\makexaxis{\thinlines
        \dcnta = \gridwidth
        \advance\dcnta by 1
        \put(\xmin, 0){\vector(1,0){\dcnta}}
        \thicklines
}
\def\makeyaxis{\thinlines
        \dcntb = \gridheight
        \advance\dcntb by 1
        \put(0, \ymin){\vector(0,1){\dcntb}}
        \thicklines
}
\def\makeaxes{\thinlines
        \dcnta = \gridwidth
        \advance\dcnta by 1
        \put(\xmin, 0){\vector(1,0){\dcnta}}
        \dcntb = \gridheight
        \advance\dcntb by 1
        \put(0, \ymin){\vector(0,1){\dcntb}}
        \thicklines
}
 
% scaling on axes
\def\makexscale{
   \dcnta = \xmin
 %\ifnum\dcnta<-9\dcnta=-9\fi
   \loop\ifnum\dcnta<0
  \put(\dcnta,0){\line(0,-1){0.1}}
      \put(\dcnta,-0.4){\makebox(0,0)[tr]{
  {\scalefont\number\dcnta}}}
      \advance\dcnta by 1
  \repeat
  \dcnta = 1
  \dcntb=\xmax
  %\ifnum\dcntb>9\dcntb=9\fi
  \loop\ifnum\dcnta<\dcntb
      \put(\dcnta,0){\line(0,-1){0.1}}
      \put(\dcnta, -0.4){\makebox(0,0)[tr]{
  {\scalefont\number\dcnta}}}
      \advance\dcnta by 1
  \repeat
  \put(\dcntb,0){\line(0,-1){0.1}}
  \put(\dcntb, -0.4){\makebox(0,0)[tr]{{\scalefont\number\dcntb}}}
}
\def\makeyscale{
   \dcnta = \ymin
   \loop\ifnum\dcnta<0
      \put(0,\dcnta){\line(-1,0){0.1}}
      \put(-0.3, \dcnta){\makebox(0,0)[r]{
  {\scalefont\number\dcnta}}}
      \advance\dcnta by 1
   \repeat
   \dcnta = 1
   \loop\ifnum\dcnta<\ymax
      \put(0,\dcnta){\line(-1,0){0.1}}
      \put(-0.3, \dcnta){\makebox(0,0)[r]{
  {\scalefont\number\dcnta}}}
      \advance\dcnta by 1
   \repeat
 \put(0,\ymax){\line(-1,0){0.1}}
   \put(-0.3, \ymax){\makebox(0,0)[r]{{\scalefont\number\ymax}}}
}

% macro for captions:
\def\gridcaption#1{
 \dcnta=\xmin
 \dcntb=\gridwidth
 \divide\dcntb by 2
 \advance\dcnta by \dcntb
 % \advance\dcnta by 1
 \dcntb=\ymax
 % \advance\dcntb by 2
 \advance\dcntb by 1
 \put(\dcnta,\dcntb){ \begin{picture}(1,1)(0,0)
        \put(0.5,0.5){\makebox(0,0){#1}}
        \end{picture}
       }
}

% macros for drawing skew dotted lines
\newsavebox{\skewdotteddownbox}
\newsavebox{\skewdottedupbox}
\savebox{\skewdotteddownbox}(1,1)[lb]{
 \put(0.2,-0.2){\circle*{0.01}}
   \put(0.4,-0.4){\circle*{0.01}}
  \put(0.6,-0.6){\circle*{0.01}}
 \put(0.8,-0.8){\circle*{0.01}}
}
\savebox{\skewdottedupbox}(1,1)[lb]{
 \put(0.2,0.2){\circle*{0.01}}
   \put(0.4,0.4){\circle*{0.01}}
  \put(0.6,0.6){\circle*{0.01}}
 \put(0.8,0.8){\circle*{0.01}}
}
\def\downskewline#1#2#3{
 \multiput(#1,#2)(1,-1){#3}{\usebox{\skewdotteddownbox}}
}
\def\upskewline#1#2#3{
 \multiput(#1,#2)(1,1){#3}{\usebox{\skewdottedupbox}}
}

% Special macros for drawing lattice paths:
\newsavebox{\hplainbox}
\savebox{\hplainbox}{\put(0,0){\line(1,0){1}}}
\newsavebox{\hdottedbox}
\savebox{\hdottedbox}(1,0)[bl]{
 \put(0.2,0){\circle*{0.01}}
 \put(0.4,0){\circle*{0.01}}
 \put(0.6,0){\circle*{0.01}}
 \put(0.8,0){\circle*{0.01}}
}
\newsavebox{\vplainupbox}
\savebox{\vplainupbox}{\put(0,0){\line(0,1){1}}}
\newsavebox{\vdottedupbox}
\savebox{\vdottedupbox}(0,1)[bl]{
 \put(0,0.2){\circle*{0.01}}
 \put(0,0.4){\circle*{0.01}}
 \put(0,0.6){\circle*{0.01}}
 \put(0,0.8){\circle*{0.01}}
}
\newsavebox{\vplaindownbox}
\savebox{\vplaindownbox}{\put(0,0){\line(0,-1){1}}}
\newsavebox{\vdotteddownbox}
\savebox{\vdotteddownbox}(0,1)[bl]{
 \put(0,-0.2){\circle*{0.01}}
 \put(0,-0.4){\circle*{0.01}}
 \put(0,-0.6){\circle*{0.01}}
 \put(0,-0.8){\circle*{0.01}}
}
\newsavebox{\skewupbox}
\savebox{\skewupbox}{\put(0,0){\line(1,1){1}}}
\newsavebox{\skewdownbox}
\savebox{\skewdownbox}{\put(0,0){\line(1,-1){1}}}
\newsavebox{\hstep}  % horizontal step
\newsavebox{\vstep}  % vertical step
\newsavebox{\sstep}  % skew step
\newsavebox{\rvstep}  % reverse vertical step
\newsavebox{\rsstep}  % reverse skew step
\newcount\updownincrement

% Macros for automatic labelling of steps:
\newif\iflabelled
\newif\ifelabel
\newcount\labelinfty
\def\setlabelinfty#1{\global\labelinfty=#1}
\newcount\labelcount
\setlabelinfty{1000}
\def\printlabel{
 \ifnum\labelcount<\labelinfty{\scalefont\number\labelcount}
 \else$\scriptstyle\infty$
 \fi
}

\def\uppath#1#2#3{
 \global\labelcount=1 % If path starts with horizontal step -> 1
 \savebox{\hstep}{\usebox{\hplainbox}}
 \savebox{\vstep}{\usebox{\vplainupbox}}
 \savebox{\sstep}{\usebox{\skewupbox}}
 \savebox{\rvstep}{\usebox{\vplaindownbox}}
 \savebox{\rsstep}{\usebox{\skewdownbox}}
 \updownincrement=1
 \dcnta=#1
 \dcntb=#2
 \thickpoint{\dcnta}{\dcntb}
 \dosteplist#3\mendlist
 \thickpoint{\dcnta}{\dcntb}
}
\def\updottedpath#1#2#3{
 \global\labelcount=1 % If path starts with horizontal step -> 1
 \savebox{\hstep}{\usebox{\hdottedbox}}
 \savebox{\vstep}{\usebox{\vdottedupbox}}
 \savebox{\sstep}{\usebox{\skewdottedupbox}}
 \savebox{\rvstep}{\usebox{\vdotteddownbox}}
 \savebox{\rsstep}{\usebox{\skewdotteddownbox}}
 \updownincrement=1
 \dcnta=#1
 \dcntb=#2
 \thickpoint{\dcnta}{\dcntb}
 \dosteplist#3\mendlist
 \thickpoint{\dcnta}{\dcntb}
}
\def\downpath#1#2#3#4{
 \global\labelcount=#4 % If path starts with horizontal step -> labeled #4
 \savebox{\hstep}{\usebox{\hplainbox}}
 \savebox{\vstep}{\usebox{\vplaindownbox}}
 \savebox{\sstep}{\usebox{\skewdownbox}}
 \savebox{\rvstep}{\usebox{\vplainupbox}}
 \savebox{\rsstep}{\usebox{\skewupbox}}
 \updownincrement=-1
 \dcnta=#1
 \dcntb=#2
 \thickpoint{\dcnta}{\dcntb}
 \dosteplist#3\mendlist
 \thickpoint{\dcnta}{\dcntb}
 \elabelfalse % for Giambelli special case only!
}
\def\downdottedpath#1#2#3#4{
 \global\labelcount=#4 % If path starts with horizontal step -> labeled #4
 \savebox{\hstep}{\usebox{\hdottedbox}}
 \savebox{\vstep}{\usebox{\vdotteddownbox}}
 \savebox{\sstep}{\usebox{\skewdotteddownbox}}
 \savebox{\rvstep}{\usebox{\vdottedupbox}}
 \savebox{\rsstep}{\usebox{\skewdottedupbox}}
 \updownincrement=-1
 \dcnta=#1
 \dcntb=#2
 \thickcircle{\dcnta}{\dcntb}
 \dosteplist#3\mendlist
 \thickcircle{\dcnta}{\dcntb}
 \elabelfalse % for Giambelli special case only!
}
% For segments of paths
\def\segment#1#2#3{
 \labelledfalse
 \savebox{\hstep}{\usebox{\hplainbox}}
 \savebox{\vstep}{\usebox{\vplainupbox}}
 \savebox{\sstep}{\usebox{\skewupbox}}
 \updownincrement=1
 \dcnta=#1
 \dcntb=#2
 \gridpoint{\dcnta}{\dcntb}
 \dosteplist#3\mendlist
 \gridpoint{\dcnta}{\dcntb}
}
\def\dottedsegment#1#2#3{
 \labelledfalse
 \savebox{\hstep}{\usebox{\hdottedbox}}
 \savebox{\vstep}{\usebox{\vdottedupbox}}
 \savebox{\sstep}{\usebox{\skewdottedupbox}}
 \updownincrement=1
 \dcnta=#1
 \dcntb=#2
 \gridcircle{\dcnta}{\dcntb}
 \dosteplist#3\mendlist
 \gridcircle{\dcnta}{\dcntb}
}
% For ``unbounded'' paths running to infinity:
\def\uparrowpath#1#2#3{
 \global\labelcount=1 % If path starts with horizontal step -> 1
 \savebox{\hstep}{\usebox{\hplainbox}}
 \savebox{\vstep}{\usebox{\vplainupbox}}
 \savebox{\sstep}{\usebox{\skewupbox}}
 \savebox{\rvstep}{\usebox{\vplaindownbox}}
 \savebox{\rsstep}{\usebox{\skewdownbox}}
 \updownincrement=1
 \dcnta=#1
 \dcntb=#2
 \thickpoint{\dcnta}{\dcntb}
 \dosteplist#3\mendlist
 \put(\dcnta,\dcntb){\vector(0,1){0.3}}
}
\def\updottedarrowpath#1#2#3{
 \global\labelcount=1 % If path starts with horizontal step -> 1
 \savebox{\hstep}{\usebox{\hdottedbox}}
 \savebox{\vstep}{\usebox{\vdottedupbox}}
 \savebox{\sstep}{\usebox{\skewdottedupbox}}
 \savebox{\rvstep}{\usebox{\vdotteddownbox}}
 \savebox{\rsstep}{\usebox{\skewdotteddownbox}}
 \updownincrement=1
 \dcnta=#1
 \dcntb=#2
 \thickcircle{\dcnta}{\dcntb}
 \dosteplist#3\mendlist
 \put(\dcnta,\dcntb){\vector(1,0){0.3}}
}

\newsavebox{\dummybox}
\def\xboxes{\savebox{\dummybox}{\usebox{\vstep}}
    \savebox{\vstep}{\usebox{\rvstep}}
    \savebox{\rvstep}{\usebox{\dummybox}}
    \savebox{\dummybox}{\usebox{\sstep}}
    \savebox{\sstep}{\usebox{\rsstep}}
    \savebox{\rsstep}{\usebox{\dummybox}}
}
    
\def\mendlist{\mendlist}
\def\dosteplist{\afterassignment\handlenextstep\let\next=}
\def\handlenextstep{
 \ifx\next\mendlist
    \let\next=\relax
 \else
  \ifx\next-
   \put(\dcnta,\dcntb){\usebox{\hstep}}
   \iflabelled\put(\dcnta,\dcntb){\makebox(1,0.6){\printlabel}}\fi
   \ifelabel\advance\labelcount by \updownincrement\fi
   \advance\dcnta by 1
  \else
   \ifx\next|
    \advance\labelcount by \updownincrement
    \put(\dcnta,\dcntb){\usebox{\vstep}}
    \advance\dcntb by\updownincrement
   \else
    \ifx\next/
     \put(\dcnta,\dcntb){\usebox{\sstep}}
     \iflabelled
      \ifnum\updownincrement>0
       \advance\labelcount by 1 % special asymmetric labelling...
       \put(\dcnta,\dcntb){\begin{picture}(1,1)(0,0)
        \put(-0.5,0.5){\printlabel}\end{picture}}
      \else
       \put(\dcnta,\dcntb){\begin{picture}(1,1)(0,0)
        \put(0.3,-0.5){\printlabel}\end{picture}}
      \fi
     \fi
     \advance\dcnta by 1
     \advance\dcntb by\updownincrement
    \else
     \ifx\next* % Up/down change for Giambelli lattice paths!
      \ifelabel\elabelfalse\else\elabeltrue\fi
      \xboxes % exchange up/down boxes, e-/h-labelling! 
      \updownincrement=-\updownincrement
     \else
      \errmessage{PATHS: Wrong symbol.}
     \fi
    \fi
   \fi
  \fi
  \let\next=\dosteplist
 \fi
 \next
}
 
% draw thick point:
\def\thickpoint#1#2{\put(#1, #2){\circle*{0.4}}}

% draw thick circle
\def\thickcircle#1#2{\put(#1, #2){\circle{0.4}}}

% draw thin circle
\def\thincircle#1#2{\put(#1, #2){\circle{0.2}}}

% draw ``medium'' point:
\def\gridpoint#1#2{\put(#1, #2){\circle*{0.2}}}
 
% draw ``medium'' circle:
\def\gridcircle#1#2{\put(#1, #2){\circle{0.2}}}

% label single horizontal step
\def\hlabel#1#2#3{\put(#1,#2){\makebox(1,0.6){\scalefont #3}}}

% label thick point above or below:
\def\labelpointal#1#2#3{       % above, to the left
 \put(#1,#2){
  \begin{picture}(1,1)(0,0)
   \put(-0.85,0.35){$\scriptstyle #3$}
  \end{picture}
 }
}

\def\labelpointar#1#2#3{       % above, to the right
 \put(#1,#2){
  \begin{picture}(1,1)(0,0)
   \put(-0.15,0.35){$\scriptstyle #3$}
  \end{picture}
 }
}

\def\labelpointbl#1#2#3{       % below, to the left
 \put(#1,#2){
  \begin{picture}(1,1)(0,0)
   \put(-0.85,-0.8){$\scriptstyle #3$}
  \end{picture}
 }
}

\def\labelpointbr#1#2#3{       % below, to the right
 \put(#1,#2){
  \begin{picture}(1,1)(0,0)
   \put(-0.15,-0.8){$\scriptstyle #3$}
  \end{picture}
 }
}

\def\endgrid{\end{picture}}

\def\drawrect#1#2#3#4{
 \put(#1,#2){%
  \begin{picture}(#3,#4)(0,0)\thicklines%
   \put(0,0){\line(0,1){#4}}
   \put(0,0){\line(1,0){#3}}
   \put(#3,#4){\line(0,-1){#4}}
   \put(#3,#4){\line(-1,0){#3}}
  \end{picture}
 }
}

\newsavebox{\smallrectbox}
\savebox{\smallrectbox}(0.2,0.2)[bl]{
 \thinlines
 \put(-0.2,-0.2){\line(1,0){0.4}}
 \put(0.2,-0.2){\line(0,1){0.4}}
 \put(0.2,0.2){\line(-1,0){0.4}}
 \put(-0.2,0.2){\line(0,-1){0.4}}
}
\def\gridrect#1#2{\put(#1, #2){\usebox{\smallrectbox}}}

%%%@@@ end input from paths

%%%@@@ input from tableaux
%%%@@@ input from tableaux
%%%%%%%%%%%%%%%%%%%%%%%%%%%%%%%%%%%%%%%%%%%%%%%%%%%%%%%%%%%%%%%%%%%%%%%%%%%%
%  TABLEAUX.TEX: Begonnen am 20. Juli 1993                                 %
% Dieses kleine Macropaket soll einfache Befehle zur Darstellung von       %
% Tableaux  einfuehren. (Nur fuer AMSLaTeX|)                               %
%%%%%%%%%%%%%%%%%%%%%%%%%%%%%%%%%%%%%%%%%%%%%%%%%%%%%%%%%%%%%%%%%%%%%%%%%%%%
\newcount\rowcount
\newcount\columncount

\def\begintableau#1#2{                           % rows, columns
\global\rowcount=#1
\begin{picture}(#2,#1)(0,0)
\thinlines
% \put(0,0){\line(0,1){#1}}
% \put(0,#1){\line(1,0){#2}}
}

\def\endtableau{\end{picture}}

\def\row#1{
   \global\advance\rowcount by -1
   \global\columncount=0
   \dorowlist#1\mendlist
}

\def\mendlist{\mendlist}

\def\dorowlist{\afterassignment\handlenextentry\let\next=}

\def\handlenextentry{
   \ifx\next\mendlist
      \let\next=\relax
   \else
      \ifx\next=-
         \put(\columncount,\rowcount){\framebox(1,1){\space}}
         \global\advance\columncount by 1
      \else
         \put(\columncount,\rowcount){\framebox(1,1){$\next$}}
         \global\advance\columncount by 1
      \fi
      \let\next=\dorowlist
   \fi
   \next
}

   % Makros zur einfachen Eingabe von Tableaux...              %
%%%@@@ end input from tableaux

%%%@@@ end input from tableaux

\section{Introduction}
\label{intro}

%%%@@@ input from parts/intro
Usually, bijective proofs of determinant identities involve
the following steps (cf., e.g, \cite[Chapter 4]{Stanton-White} or \cite{Zeilberger,Zeilberger2}):
\begin{itemize}
\item Expansion of the determinant as sum over the symmetric group,
\item Interpretation of this sum as the generating function of some
        set of combinatorial objects which are equipped with some signed weight,
\item Construction of an explicit weight-- and sign--preserving bijection
        between the respective combinatorial objects, maybe supported by the
        construction of a sign--reversing involution for certain objects.
\end{itemize}

Here, we will present another ``method'' of bijective proofs for
determinant identitities, which involves the following steps:
\begin{itemize}
\item First, we replace the entries $a_{i,j}$ of the determinants by
        $h_{\lambda_i-i+j}$ (where $h_m$ denotes the $m$--th complete
        homogeneous function),
\item Second, by the Jacobi--Trudi identity we transform the original
        determinant identity into an equivalent identity for Schur functions,
\item Third, we obtain a bijective proof for this equivalent
        identity by using the interpretation of Schur functions in
		  terms of nonintersecting lattice paths.
		  (In this paper, we shall achieve this with a
        construction which was used for the proof of a Schur function
        identity \cite[Theorem~1.1]{Fulmek:Ciucu} conjectured by Ciucu.)
\end{itemize}

We show how this method applies naturally to provide elegant bijective
proofs of Dodgson's Condensation Rule \cite{Dodgson}
and the Pl\"ucker relations.

The bijective construction we use here was (to the best of our knowledge) first used by I.~Goulden \cite{Goulden:Schur}.
(The first author is grateful to A.~Hamel \cite{Hamel}
for drawing his attention to Goulden's work.)
Goulden's exposition, however, left open a small gap, which we shall close
here.

The paper is organized as follows: In Section~\ref{expo}, we
present the theorems we want to prove, and explain Steps 1 and 2
of our above ``method'' in greater detail. In
Section~\ref{background}, we briefly recall the combinatorial
definition of Schur functions and the Gessel--Viennot--approach.
In Section~\ref{bijection}, we explain the bijective construction
employed in Step 3 of our ``method'' by
using the proof of a Theorem from Section~\ref{expo} as an illustrating
example. There, we shall also close the small gap in Goulden's
work. In Section~\ref{general}, we ``extract'' the general
structure underlying the bijection: As it turns out, this is just a
simple graph--theoretic statement. From this we may easily derive a
general ``class'' of Schur function identities which follow from
these considerations. In order to show that these quite general
identitities specialize to something useful, we shall deduce the
Pl\"ucker relations, using again our ``method''. In
Section~\ref{kleber}, we turn to a theorem
% markus
\cite[Theorem 3.2]{Kleber}
recently proved by the second author by using Pl\"ucker relations: We explain how this theorem fits into our
construction and give a bijective proof using inclusion--exclusion.

%%%@@@ end input from parts/intro

\section{Exposition of identities and proofs}
\label{expo}

%%%@@@ input from parts/exposition
The origin of this paper was the attempt to give a bijective proof of
the following identity for Schur functions, which arose in work of Kirillov \cite{Kirillov}:

\begin{thm}
\label{thm:kirillov}
Let $c,r$ be positive integers; denote by $\qp{c}{r}$ the partition
consisting of $r$
rows with constant length $c$. Then we have the following identity
for Schur functions:
\begin{equation}
\label{eq:kirillov}
\left(\sqp{c}{r}\right)^2=
        \sqp{c}{r-1}\cdot\sqp{c}{r+1}+\sqp{(c-1)}{r}\cdot\sqp{(c+1)}{r}.
\end{equation}
\end{thm}

(See \cite[7.10]{Stanley2}, \cite{Fulton-Harris:Representation-Theory}, \cite{Macdonald:Symmetric-Functions} or \cite{Sagan:The-Symmetric}
for background information on Schur functions; in order to keep our exposition
self--contained, a combinatorial definition is given in
Section~\ref{bijection}.)

The identity \eqref{eq:kirillov} was recently considered by
% markus
the second author
\cite[Theorem 4.2]{Kleber}, who also gave a bijective proof,
and generalized it considerably \cite[Theorem 3.2]{Kleber}.

The construction we use here
does in fact prove a more general statement:

\begin{thm}
\label{thm:general}
Let $\pa{\lambda_1,\lambda_2,\dots,\lambda_{r+1}}$ be a
partition, where $r>0$ is some integer.
Then we have the following identity
for Schur functions:
\begin{multline}
\label{eq:general}
\sla{\lambda_1,\dots,\lambda_r}\cdot\sla{\lambda_2,\dots,\lambda_{r+1}}\\
        =\sla{\lambda_2,\dots,\lambda_r}\cdot\sla{\lambda_1,\dots,\lambda_{r+1}} +
        \sla{\lambda_2-1,\dots,\lambda_{r+1}-1}\cdot
        \sla{\lambda_1+1,\dots,\lambda_r+1}.
\end{multline}
\end{thm}

Clearly, Theorem~\ref{thm:kirillov} is a direct consequence of Theorem~\ref{thm:general}: Simply set $\lambda_1=\dots=\lambda_{r+1}=c$.

Theorem~\ref{thm:general}, however, is in fact equivalent to Dodgson's
condensation formula \cite{Dodgson}, which is also known as
Desnanot--Jacobi's adjoint matrix theorem (see \cite[Theorem 3.12]{Bressoud}:
According to \cite{Bressoud}, Langrange discovered this theorem for $n=3$,
% Lagrange: 1773
Desnanot proved it for $n\leq 6$
% Desnanot: 1819
and Jacobi published the general theorem \cite{Jacobi},
% Jacobi: in 1833, paper 1841
% but first proved in 1819 by P.~Desnanot, according to
see also \cite[vol.~I, pp.~142]{MuirAB}):

\begin{thm}
\label{lem:general-minors}
Let $A$ be an arbitrary $(r+1)\times(r+1)$--determinant. Denote by
        $\minor{A}{r_1}{r_2}{c_1}{c_2}$
the minor consisting of rows $r_1,r_1+1,\dots,r_2$ and columns
$c_1,c_1+1,\dots,c_2$ of $A$. Then we have the following identity:
\begin{multline}
\label{eq:general-minors}
\minor{A}{1}{r+1}{1}{r+1}\minor{A}{2}{r}{2}{r} \\=
\minor{A}{1}{r}{1}{r}\minor{A}{2}{r+1}{2}{r+1} -
\minor{A}{2}{r+1}{1}{r}\minor{A}{1}{r}{2}{r+1}.
\end{multline}
\end{thm}

The transition from Theorem~\ref{lem:general-minors} to
Theorem~\ref{thm:general} is established by the Jacobi--Trudi identity (see
\cite[I, (3.4)]{Macdonald:Symmetric-Functions}), which
states that for any partition $\lambda=(\lambda_1,\dots,\lambda_r)$
of length $r$ we have
\begin{equation}
\label{eq:Jacobi-Trudi}
s_\lambda = \det(h_{\lambda_i-i+j})_{i,j=1}^r\;,
\end{equation}
where $h_m$ denotes the $m$--th complete homogeneous symmetric function:
% \begin{equation}
% h=
% \end{equation}
Setting $A_{i,j}\defeq h_{\lambda_i-i+j}$ for $1\leq i,j\leq r+1$ in
Theorem~\ref{lem:general-minors} and using identity \eqref{eq:Jacobi-Trudi}
immediately yields \eqref{eq:general}.

That the seemingly weaker statement of Theorem~\ref{thm:general} does
in fact imply Theorem~\ref{lem:general-minors} is due to the
following observation: Choose $\lambda$ so that the numbers
$\lambda_i-i+j$ are all distinct for $1\leq i,j\leq (r+1)$ (e.g.,
$\lambda=\left((r+1)r,r^2,(r-1)r,\dots,r\right)$ would suffice)
% in Theorem~\ref{thm:general}
and rewrite \eqref{eq:general} as a
determinantal expression according to the Jacobi--Trudi identity
\eqref{eq:Jacobi-Trudi}.
This yields a special case of identity
\eqref{eq:general-minors}
with $A_{i,j}\defeq h_{\lambda_i-i+j}$ as above.
% Note that all the
% entries $A_{i,j}$ are distinct for this special choice of
% $\lambda$, and recall
Now recall that the complete homogeneous symmetric functions are
algebraically independent (see, e.g., \cite{Sturmfels}), whence the
identity
\eqref{eq:general-minors} is true for generic $A_{i,j}$. For later
use, we record this simple observation in a more general fashion:
\begin{obs}
\label{obs:christian}
Let ${\mathcal I}$ be an identity involving determinants
of homogeneous symmetric
functions $h_{n}$, where $n$ is some nonnegative integer.
% Suppose that
% ${\mathcal I}$ can be specialized so that each of the indices $n$
% appears only once in any determinant involved:
Then ${\mathcal I}$ is, in fact, equivalent to a
general determinant identity which is obtained from ${\mathcal I}$ by considering
each $h_n$ as a formal variable.
\end{obs}

So far, the promised proof (to be given in Section~\ref{bijection})
of Theorem~\ref{thm:general}
would give a new bijective proof of Dodgson's
Determinant--Evaluation Rule (a beautiful bijective proof was also
given by Zeilberger \cite{Zeilberger}). But we can do a little
better: Our bijective construction does, in fact, apply to a quite
general ``class of Schur function identities'', a special case of
which implies the Pl\"ucker relations (also known as
Grassmann--Pl\"ucker syzygies), see, e.g., \cite{Sturmfels}, or
\cite[Chapter 3, Section 9, formula II]{Turnbull}:
\begin{thm}[Pl\"ucker relations]
\label{thm:pluecker}
Consider an arbitrary $2n\times n$--matrix with row indices
$1,2,\dots,2n$.
Denote the $n\times n$--minor
of this matrix consisting of rows $i_1,\dots,i_n$ by $[i_1,\dots,i_n]$.

Consider some fixed
list of integers $1\leq r_1< r_2<\dots< r_k\leq n$, $0\leq k\leq n$.
Then we have:
\begin{multline}
\label{eq:pluecker}
[1,2,\dots,n]\cdot [n+1,n+2,\dots,2n]=\\
        \sum_{n+1\leq s_1< s_2<\dots< s_k\leq 2n}
                [1,\dots,s_1,\dots,s_k,\dots,n]\cdot
                [n+1,\dots,r_1,\dots,r_k,\dots,2n],
\end{multline}
where the notation of the summands means that rows $r_i$ were exchanged
with rows $s_i$, respectively.
\end{thm}

This is achieved by observing that \eqref{eq:pluecker} can be specialized
to a Schur function identity of the form
\begin{equation*}
\label{eq:pluecker-schur}
s_\lambda s_\mu=
        \sum_{\lambda^\prime\!\!,\,\mu^\prime} s_{\lambda^\prime}s_{\mu^\prime},
\end{equation*}
where $\lambda$ and $\mu$ are partitions with the same number $n$ of parts,
and where the sum is over certain pairs $\lambda^\prime,\mu^\prime$ derived
from $\lambda, \mu$ (to be described later). This Schur function identity
belongs to the ``class of identities'' which follow from the bijective
construction. By applying Observation~\ref{obs:christian} with
suitable $\lambda$ and $\mu$, we may deduce \eqref{eq:pluecker}.

\begin{rem}
Summing equation~\eqref{eq:pluecker} over all possible choices of
subsets $\{r_1,\dots,r_k\}$ yields the determinant identity behind
Ciucu's Schur function identity \cite[Theorem~1.1]{Fulmek:Ciucu}
\begin{equation}
\label{eq:ciucu}
\sum_{A\subset T:\;\vert A\vert =k} s_{\lambda(A)}s_{\lambda(T-A)} =
	2^k s_{\lambda(t_2,\dots,t_{2k})}s_{\lambda(t_1,\dots,t_{2k-1})},
\end{equation}
where $T=\{t_1<\dots<t_{2k}\}$ is some set of positive integers
and $\lambda(\{t_{i_1}<\dots <t_{i_r}\})$ denotes the partition with parts
$t_{i_r}-r+1\geq\dots\geq t_{i_2}-1\geq t_{i_1}$.
\end{rem}

\begin{rem}
The Pl\"ucker relations \eqref{eq:pluecker} appear in a slightly different
notation as Theorem~2 in \cite{Stanley-Propp}, together with another
elegant proof.
\end{rem}

% markus
Moreover, the bijective method yields a proof of the second author's
theorem \cite[Theorem 3.2]{Kleber}: Since this theorem is rather
complicated to state, we defer it to Section~\ref{kleber}.

%%%@@@ end input from parts/exposition

\section{Combinatorial background and definitions}
\label{background}

%%%@@@ input from parts/background
As usual,
an $r$-tuple $\lambda = \parof{\lambda_1, \lambda_2,\dots,
\lambda_r}$ with $\lambda_1\geq\lambda_2\geq\dots\geq\lambda_r\geq 0$
is called a {\em partition of length $r$\/}. The {\em Ferrers
board\/} $F(\lambda)$ of $\lambda$
is an array of cells with $r$ left-justified rows and $\lambda_i$
cells in row $i$.

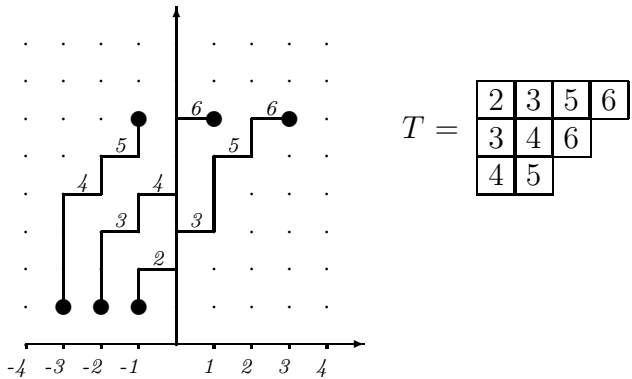
\begin{figure}
\caption{Illustration of a $6$--semistandard Young tableau and
its associated lattice paths for $\lambda=(4,3,2)$.}
\label{fig:GV}
\begin{picture}(25,12)(-4,0)
\put(12,5){
\begintableau{3}{4}
\row{2356}
\row{346}
\row{45}
\put(-2,1.5){$T$ = }
\endtableau
}
\put(0,1){
\begingrid{-4}{0}{4}{8}
\makeaxes
\makexscale

\elabelfalse
\labelledtrue

% Die Pfade:
\uppath{-1}{1}{|-|-||-|-}
\uppath{-2}{1}{||-|-||-}
\uppath{-3}{1}{|||-|-|}
\endgrid
}
\end{picture}
\end{figure}

An {\em $N$--semistandard Young tableau\/} of shape $\lambda$ is a filling of the cells of $F(\lambda)$ with integers from the set $\{1,2,\dots,N\}$, such
that the numbers filled into the cells weakly increase in rows
and strictly increase in columns
(see the right picture of Figure~\ref{fig:GV} for an illustration).

Schur functions, which are irreducible general linear characters, can be combinatorially defined by means of $N$--semistandard Young tableaux (see
\cite[I, (5.12)]{Macdonald:Symmetric-Functions},
\cite[Def.~4.4.1]{Sagan:The-Symmetric}, \cite[Def.~5.1]{Stanley}):

\begin{equation*}
s_\lambda(x_1, x_2, x_3, \dots, x_N) = \sum_{{\bold T}}w({\bold T}),
\end{equation*}
where the sum is over all $N$--semistandard Young tableaux ${\bold T}$ of
shape $\lambda$. Let $m({\bold T}, k)$ be the number of entries $k$
in the tableau ${\bold T}$. The weight $w({\bold T})$ of ${\bold T}$
is defined as follows:
\begin{equation*}
w({\bold T}) = \prod_{k = 1}^{N} x_k^{m({\bold T}, k)}.
\end{equation*}

The Gessel-Viennot interpretation \cite{Gessel-Viennot:Determinants-paths}
of semistandard Young tableaux of shape $\lambda$ as nonintersecting lattice paths (see the left picture of Figure~\ref{fig:GV} for an illustration)
allows an equivalent definition of Schur functions:

\begin{equation*}
s_\lambda(x_1, x_2, x_3, \dots, x_N) = \sum_{{\bold P}}w({\bold P}),
\end{equation*}
where the sum is over all
$r$-tuples ${\bold P} = \left(P_1, P_2, \dots, P_r\right)$ of
lattice paths (in the integer lattice, i.e., the directed graph with
vertices $\mathbb Z\times\mathbb Z$ and arcs from $(j, k)$ to $(j+1, k)$ and
from $(j, k)$ to $(j, k+1)$ for all $j, k$), where $P_i$ starts at $(-i, 1)$
and ends at $(\lambda_{i}-i, N)$, and where no two paths $P_i$ and $P_j$
have a lattice point in common (such an $r$-tuple is called nonintersecting).

The weight $w({\bold P})$ of
an % nonintersecting
$r$-tuple ${\bold P} = \left(P_1, P_2, \dots, P_r\right)$ of paths
is defined by:
\begin{equation*}
w({\bold P}) = \prod_{i = 1}^{r} w({P_i}).
\end{equation*}

The weight $w(P)$ of a single path $P$ is defined as follows:
Let $n({P}, k)$ be the number of horizontal steps at height $k$
(i.e., directed arcs from some $(j, k)$ to $(j+1, k)$) that belong to
path $P$, then we define
\begin{equation*}
w({P}) = \prod_{k = 1}^{N} x_k^{n({P}, k)}.
\end{equation*}

That these definitions are in fact equivalent is due to a weight--preserving
bijection between tableaux and nonintersecting lattice paths. The
Gessel--Viennot method \cite{Gessel-Viennot:Determinants-paths}
builds on the lattice path definition to give a bijective proof of
the Jacobi--Trudi identity \eqref{eq:Jacobi-Trudi}  (see, e.g.,
\cite[ch.~4]{Sagan:The-Symmetric}, \cite{Stembridge} or
\cite{Fulmek-Krattenthaler}).

Next, we give a combinatorial definition for {\em skew\/} Schur
functions: Let $\lambda=(\lambda_1,\dots,\lambda_r)$ and
$\mu=(\mu_1,\dots,\mu_r)$ be partitions with $\mu_i\leq\lambda_i$
for $1\leq i\leq r$; here, we allow $\mu_i=0$.

The {\em skew Ferrers
board\/} $F(\lambda/\mu)$ of $(\lambda,\mu)$
is an array of cells with $r$ left-justified rows and $\lambda_i-\mu_i$
cells in row $i$, where the first $\mu_i$ cells in row $i$ are missing.

\begin{figure}
\caption{Illustration of a $6$--semistandard skew Young tableau and
its associated lattice paths
for $\lambda=(4,3,2)$ and $\mu=(1,0,0)$.}
\label{fig:skew}
\begin{picture}(25,12)(-4,0)
\put(12,5){
\begintableau{3}{4}
\row{-356}
\row{346}
\row{45}
\put(-2,1.5){$T$ = }
\endtableau
}
\put(0,1){
\begingrid{-4}{0}{4}{8}
\makeaxes
\makexscale

\elabelfalse
\labelledtrue

% Die Pfade:
\uppath{-0}{1}{||-||-|-}
\uppath{-2}{1}{||-|-||-}
\uppath{-3}{1}{|||-|-|}
\endgrid
}
\end{picture}
\end{figure}
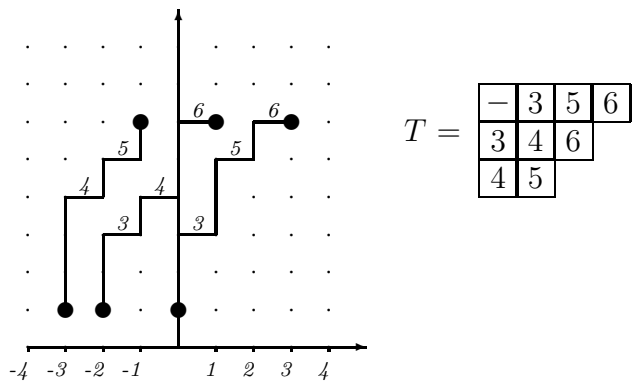

An {\em $N$--semistandard skew Young tableau\/} of shape $\lambda/\mu$ is a filling of the cells of $F(\lambda/\mu)$ with integers from the set
$\{1,2,\dots,N\}$,
such that the numbers filled into the cells weakly increase in rows and
strictly increase in columns
(see the right picture of Figure~\ref{fig:skew} for an illustration).

Then we have the following definition for skew Schur functions:

\begin{equation*}
s_{\lambda/\mu}(x_1, x_2, x_3, \dots, x_N) = \sum_{{\bold T}}w({\bold T}),
\end{equation*}
where the sum is over all $N$--semistandard skew Young tableaux ${\bold T}$ of
shape $\lambda/\mu$, where the weight $w({\bold T})$ of ${\bold T}$
is defined as before.

Equivalently, we may define:
\begin{equation*}
s_{\lambda/\mu}(x_1, x_2, x_3, \dots, x_N) = \sum_{{\bold P}}w({\bold P}),
\end{equation*}
where the sum is over all $r$-tuples ${\bold P} = \left(P_1, P_2, \dots,
P_r\right)$ of nonintersecting lattice paths, where $P_i$ starts at
$(\mu_i-i, 1)$ and ends at $(\lambda_{i}-i, N)$ (see the left
picture of Figure~\ref{fig:skew} for an illustration), and where the
weight $w({\bold P})$ of such an $r$-tuple ${\bold P}$ is defined as before.

%%%@@@ end input from parts/background

\section{Bijective proof of Theorem~\ref{thm:general}}
\label{bijection}

%%%@@@ input from parts/bijection
\begin{proof}
Let us start with a combinatorial description for the objects
involved in \eqref{eq:general}: By the Gessel--Viennot
interpretation of Schur functions as generating functions of
nonintersecting lattice paths, we may view the left--hand side of
the equation as the weight of all {\em pairs\/} $({\bold\TP}^g,
{\bold\TP}^b)$, where ${\bold\TP}^g$ and ${\bold\TP}^b$ are
$r$-tuples of nonintersecting lattice paths. The paths of
${\bold\TP}^g$ are coloured green, the paths of ${\bold\TP}^b$ are
coloured blue. The $i$-th green path $\TP^g_i$ starts at $(-i, 1)$
and ends in $(\lambda_{i}-i, N)$. The $i$-th blue path $\TP^b_i$
starts at $(-i-1, 1)$ and ends in $(\lambda_{i+1}-i-1, N)$. For an
illustration, see the upper left pictures in
Figures~\ref{fig:caseA} and \ref{fig:caseB}, where green paths are
drawn with full lines and blue paths are drawn with dotted lines.

For the right--hand side of \eqref{eq:general}, we use the same
interpretation. We may view the first term as the weight of all
{\em pairs\/} $({\bold\TQ}^g, {\bold\TQ}^b)$, where ${\bold\TQ}^g$
is an $(r-1)$-tuple of nonintersecting lattice paths and
${\bold\TQ}^b$ is an $(r+1)$-tuple of nonintersecting lattice
paths. The paths of ${\bold\TQ}^g$ are coloured green, the paths of
${\bold\TQ}^b$ are coloured blue. The $i$-th green path $\TQ^g_i$
starts at $(-i-1, 1)$ and ends in $(\lambda_{i+1}-i-1, N)$. The
$i$-th blue path $\TQ^b_i$ starts at $(-i, 1)$ and ends in
$(\lambda_{i}-i, N)$. For an illustration, see the upper right
picture in Figure~\ref{fig:caseA}.

In the same way, we may view the second term as the weight of all
{\em pairs\/} $({\bold\TS}^g, {\bold\TS}^b)$, where ${\bold\TS}^g$
and ${\bold\TS}^b$ are $r$-tuples of nonintersecting lattice paths.
The paths of ${\bold\TS}^g$ are coloured green, the paths of
${\bold\TS}^b$ are coloured blue. The $i$-th green path $\TS^g_i$
starts at $(-i, 1)$ and ends in $(\lambda_{i+1}-i-1, N)$. The
$i$-th blue path $\TS^b_i$ starts at $(-i-1, 1)$ and ends in
$(\lambda_{i}-i, N)$. For an illustration, see the upper right
picture in Figure~\ref{fig:caseB}.

In any case, the weight of some pair of paths $({\bold P},{\bold Q})$ is defined as follows:
\begin{equation*}
w({\bold P},{\bold Q}):=
        w({\bold P})\cdot w({\bold Q}).
\end{equation*}

\begin{figure}
\caption{Illustration of the construction in the proof, case A: $r=3$,
$(\lambda_1,\lambda_2,\lambda_3,\lambda_4)=(5,4,3,2)$.
% The
% upper left picture shows some
% pair $({\bold\TP}^g,{\bold\TP}^b)$, and the upper right picture shows the
% corresponding pair $({\bold\TQ}^g,{\bold\TQ}^b)$.
%
% In both pictures, blue paths start in $(i,1)$ and end in
% $(\lambda_{r+2-i}+i,N)$, green paths start in $(i+1,1)$ and end in
% $(\lambda_{r+1-i}+i+1,N)$.
%
% The orientation
% of the graph is indicated by small arrows. The ``recolouring path''
% used in the bijection is shown below: It connects the rightmost
% endpoint $(9,5)$ with the rightmost starting point $(4,1)$ of the lattice
% paths.
}
\label{fig:caseA}
 \begin{picture}(30,20)(0,0)
  \put(1,11){
   \begingrid{-5}{0}{5}{6}
    \gridcaption{Some pair $({\bold\TP}^g,{\bold\TP}^b)$:}
    \makeaxes
    \makexscale
    \makeyscale
    % The green paths:
    \uppath{-3}{1}{|||-|--}
    \uppath{-2}{1}{|-|-|-|-}
    \uppath{-1}{1}{----||-||}
    % The blue paths:
    \updottedpath{-4}{1}{||||--}
    \updottedpath{-3}{1}{|-||--|}
    \updottedpath{-2}{1}{-|--||-|}
         % The orientation of edges:
    \thinlines
    \put(4.2,4.8){\vector(0,-1){0.6}}
    \put(2.2,4.2){\vector(0,1){0.6}}
    \put(1.8,5.2){\vector(-1,0){0.6}}
    \put(0.2,4.2){\vector(0,1){0.6}}
    \put(-0.2,5.2){\vector(-1,0){0.6}}
    \put(-2.8,5.2){\vector(1,0){0.6}}
   \endgrid
  }
  \put(13,11){
   \begingrid{-5}{0}{5}{6}
    \gridcaption{Corresponding pair $({\bold\TQ}^g,{\bold\TQ}^b)$:}
    \makeaxes
    \makexscale
    \makeyscale
    % The green paths:
    \uppath{-3}{1}{|-|||--}
    \uppath{-2}{1}{-||-|-|-}
    % The blue paths:
    \updottedpath{-4}{1}{||||--}
    \updottedpath{-3}{1}{|||---|}
    \updottedpath{-2}{1}{|---||-|}
    \updottedpath{-1}{1}{----||-||}
         % The orientation of edges:
    \thinlines
    \put(4.2,4.8){\vector(0,-1){0.6}}
    \put(2.2,4.8){\vector(0,-1){0.6}}
    \put(1.2,5.2){\vector(1,0){0.6}}
    \put(0.2,4.8){\vector(0,-1){0.6}}
    \put(-0.8,5.2){\vector(1,0){0.6}}
    \put(-2.2,5.2){\vector(-1,0){0.6}}
   \endgrid
  }
  \put(7,1){
   \begingrid{-5}{0}{5}{6}
    \gridcaption{The changing trail starting in $(4,5)$:}
    \makeaxes
    \makexscale
    \makeyscale
    % The recoulouring path 1:
    \thickcircle{4}{5}
         \put(4,4.8){\line(0,-1){1.8}}
         \put(4,3){\line(-1,0){1}}
         \put(3,3){\line(0,-1){2}}
         \put(3,1){\line(-1,0){4}}
         \put(-1,1){\line(0,1){1}}
         \put(-1,2){\line(-1,0){0.7}}
         \put(-1.7,2.3){\oval(0.6,0.6)[lb]}
         \put(-2,2.3){\line(0,1){1.7}}
         \put(-2,4){\line(-1,0){1}}
         \put(-3,4){\line(0,-1){2}}
         \put(-3,2){\line(1,0){0.7}}
         \put(-2.3,1.7){\oval(0.6,0.6)[rt]}
         \put(-2,1.7){\line(0,-1){0.7}}
         \put(-2,1){\vector(1,0){0.8}}
    \thickcircle{-1}{1}
   \endgrid
  }
 \end{picture}
\end{figure}

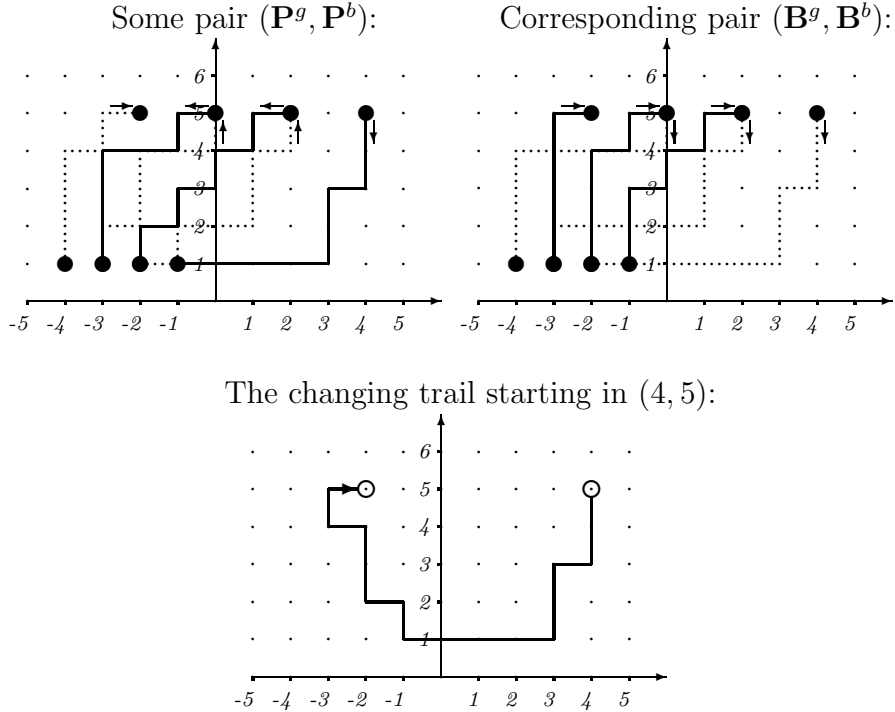
\begin{figure}
\caption{Illustration of the construction in the proof, case B: $r=3$,
$(\lambda_1,\lambda_2,\lambda_3,\lambda_4)=(5,4,3,2)$.
% ; green paths
% are drawn with full lines, blue paths are drawn with dotted lines.
% The upper left picture shows some pair
% $({\bold\TP}^g,{\bold\TP}^b)$, and the upper right picture shows
% the corresponding pair $({\bold\TS}^g,{\bold\TS}^b)$.
%
% In the left picture, blue paths start in $(i,1)$ and end in
% $(\lambda_{r+2-i}+i,N)$, green paths start in $(i+1,1)$ and end in
% $(\lambda_{r+1-i}+i+1,N)$.
% In the right picture, blue paths start in $(i,1)$ and end in
% $(\lambda_{r+1-i}+i+1,N)$, green paths start in $(i+1,1)$ and end in
% $(\lambda_{r+2-i}+i,N)$.
%
% The orientation
% of the graph is indicated by small arrows. The ``recolouring path''
% used in the bijection is shown below: It connects the rightmost
% endpoint $(9,5)$ with the leftmost endpoint $(3,5)$ of the lattice
% paths.
}
\label{fig:caseB}
 \begin{picture}(30,20)(0,0)
  \put(1,11){
   \begingrid{-5}{0}{5}{6}
    \gridcaption{Some pair $({\bold\TP}^g,{\bold\TP}^b)$:}
    \makeaxes
    \makexscale
    \makeyscale
    % The green paths:
    \uppath{-3}{1}{|||--|-}
    \uppath{-2}{1}{|-|-|-|-}
    \uppath{-1}{1}{----||-||}
    % The blue paths:
    \updottedpath{-4}{1}{|||-|-}
    \updottedpath{-3}{1}{|-||--|}
    \updottedpath{-2}{1}{-|--||-|}
         % The orientation of edges:
    \thinlines
    \put(4.2,4.8){\vector(0,-1){0.6}}
    \put(2.2,4.2){\vector(0,1){0.6}}
    \put(1.8,5.2){\vector(-1,0){0.6}}
    \put(0.2,4.2){\vector(0,1){0.6}}
    \put(-0.2,5.2){\vector(-1,0){0.6}}
    \put(-2.8,5.2){\vector(1,0){0.6}}
   \endgrid
  }
  \put(13,11){
   \begingrid{-5}{0}{5}{6}
    \gridcaption{Corresponding  pair $({\bold\TS}^g,{\bold\TS}^b)$:}
    \makeaxes
    \makexscale
    \makeyscale
    % The green paths:
    \uppath{-3}{1}{||||-}
    \uppath{-2}{1}{|||-|-}
    \uppath{-1}{1}{||-|-|-}
    % The blue paths:
    \updottedpath{-4}{1}{|||----|}
    \updottedpath{-3}{1}{|----||-|}
    \updottedpath{-2}{1}{-----||-||}
         % The orientation of edges:
    \thinlines
    \put(4.2,4.8){\vector(0,-1){0.6}}
    \put(2.2,4.8){\vector(0,-1){0.6}}
    \put(1.2,5.2){\vector(1,0){0.6}}
    \put(0.2,4.8){\vector(0,-1){0.6}}
    \put(-0.8,5.2){\vector(1,0){0.6}}
    \put(-2.8,5.2){\vector(1,0){0.6}}
   \endgrid
  }
  \put(7,1){
   \begingrid{-5}{0}{5}{6}
    \gridcaption{The changing trail starting in $(4,5)$:}
    \makeaxes
    \makexscale
    \makeyscale
    % The recoulouring path 1:
    \thickcircle{4}{5}
         \put(4,4.8){\line(0,-1){1.8}}
         \put(4,3){\line(-1,0){1}}
         \put(3,3){\line(0,-1){2}}
         \put(3,1){\line(-1,0){4.0}}
         \put(-1,1){\line(0,1){1.0}}
         \put(-1,2){\line(-1,0){1.0}}
         %
         % \put(-1.7,2.3){\oval(0.6,0.6)[lb]}
         \put(-2,2){\line(0,1){2.0}}
         \put(-2,4){\line(-1,0){1.0}}
         %
         % \put(-2.7,4.3){\oval(0.6,0.6)[lb]}
         \put(-3,4){\line(0,1){1.0}}
         \put(-3,5){\vector(1,0){0.8}}
    \thickcircle{-2}{5}
   \endgrid
  }
 \end{picture}
\end{figure}

What we want to do is to give a weight--preserving bijection between the
objects on the left side and on the right side:
\begin{equation}
\label{eq:bijection}
\{({\bold\TP}^g, {\bold\TP}^b)\}\leftrightarrow
        \left(\{({\bold\TQ}^g, {\bold\TQ}^b)\}\cup
        \{({\bold\TS}^g, {\bold\TS}^b)\}\right).
\end{equation}
Clearly, such a bijection would establish \eqref{eq:general}.

The basic idea is very simple and was already used in
\cite{Goulden:Schur} and in \cite{Fulmek:Ciucu}: Since it will be
reused later, we state it here quite generally:
\begin{dfn}
\label{dfn:graph-auxil}
Let ${\bold\TP}^1, {\bold\TP}^2$ be two arbitrary families of
nonintersecting lattice paths. The paths $\TP^1_i$ of the first
family are coloured with colour blue, the paths $\TP^2_j$ of the
second familiy are coloured with colour green.

Let $G({\bold\TP}^1, {\bold\TP}^2)$ be the ``two--coloured'' graph
made up by ${\bold\TP}^1$ and ${\bold\TP}^2$ in the obvious sense.
Observe that there are the two possible orientations for any edge
in that graph: When traversing some path, we may either move
``right--upwards'' (this is the ``original'' orientation of the
paths) or ``left--downwards''.

A {\em changing trail} is a trail in $G({\bold\TP}^1, {\bold\TP}^2)$
with the following properties:
\begin{itemize}
\item Subsequent edges of the same colour are traversed in the same
		orientation, subsequent edges of the opposite colour are traversed
		in the opposite orientation.
\item At every intersection of green and blue paths,
		colour {\em and\/} orientation are changed {\em if this is possible\/}
		(i.e., if there is an adjacent edge of opposite colour and
		opposite orientation);
		otherwise the trail must stop there.
\item The trail is {\em maximal\/} in the sense that
		it cannot be extended by adjoining edges (in a way which is
		consistent with the above conditions) at its
		start or end.
\end{itemize}

Note that for every edge $e$, there is a {\em unique\/} changing trail which
contains $e$: E.g., consider some blue edge which is right-- or
upwards--directed and enters vertex $v$. If there is an intersection at $v$,
and if there is a green edge leaving $v$ (in opposite direction left or
downwards), then the trail must continue with this edge; otherwise it must stop
at $v$. If there is no intersection at $v$,
and if there is a blue edge leaving $v$ (in the same direction right or
upwards), then the trail must continue with this edge; otherwise it must stop
at $v$.

Note that a changing trail is either
``path--like'', i.e., has obvious starting point and end point
% markus
(clearly, these must be the end points or starting points of some path from
either ${\bold\TP}^1$ or ${\bold\TP}^2$), 
or it is ``cycle--like'', i.e., is a closed trail.
\end{dfn}

Let us return from general definitions to our concrete case: Starting
with an object $({\bold\TP}^g, {\bold\TP}^b)$ from the left--hand side of
\eqref{eq:bijection}, we interpret this pair of lattice paths
as a graph $G({\bold\TP}^g, {\bold\TP}^b)$ with
green and blue edges. (See the upper left
pictures in Figures~\ref{fig:caseA} and \ref{fig:caseB}.)

Next, we determine the changing trail which starts at the rightmost
endpoint $(\lambda_{1}-1,N)$: Follow the green edges downward or to
the left; at every intersection, change colour and orientation, if this
is possible; otherwise stop there.
Clearly, this changing trail is ``path--like''.
(See Figures~\ref{fig:caseA} and
\ref{fig:caseB} for an illustration: There, the orientation of edges
is indicated by small arrows in the upper pictures; the lower pictures
show the corresponding changing trails.)

Now we change colours green to blue and vice versa along this
changing trail: It is easy to see that this recolouring
yields nonintersecting tuples of green and blue lattice paths.

Note that there are exactly two possible cases:

\noindent{\bf Case A:} The changing trail stops at the rightmost
starting point, $(-1,1)$, of the lattice paths. In this case, from
the recolouring procedure we obtain an object $({\bold\TQ}^g,
{\bold\TQ}^b)$; see the upper right picture in
Figure~\ref{fig:caseA}.

\noindent{\bf Case B:} The changing trail stops at the the leftmost endpoint,
$(\lambda_{r+1}-r-1,N)$, of the lattice paths. In this case, from
the recolouring procedure we obtain an object $({\bold\TS}^g,
{\bold\TS}^b)$; see the upper right picture in
Figure~\ref{fig:caseB}.

It is clear that altogether this gives a mapping of the set of all
objects $({\bold\TP}^g, {\bold\TP}^b)$ into the union of the two sets of all
objects $({\bold\TQ}^g, {\bold\TQ}^b)$ and $({\bold\TS}^g, {\bold\TS}^b)$,
respectively. Of course, this mapping is weight--preserving. It is also
injective since the above construction is reversed by simply repeating
it, i.e, determine the changing trail starting  at the rightmost
endpoint $(\lambda_{1}-1,N)$ (this trail is exactly the same as before,
only the colours are exchanged) and change colours.
For an illustration, read Figures~\ref{fig:caseA} and
\ref{fig:caseB} from right to left.

So what is left to prove is surjectivity: To this end, it suffices
to prove that if we apply our (injective) recolouring construction
to an {\em arbitrary\/} object $({\bold\TQ}^g, {\bold\TQ}^b)$ or
$({\bold\TS}^g, {\bold\TS}^b)$, we do {\em always\/} get an object
$({\bold\TP}^g, {\bold\TP}^b)$; i.e., two r--tuples of
nonintersecting lattice paths, coloured green and blue, and with
the appropriate starting points and endpoints.

We {\em do\/} have something to prove: Note that in both cases, A
(see Figure~\ref{fig:caseA}) and B (see Figure~\ref{fig:caseB}),
there is {\em prima vista\/} a second possible endpoint for the
changing trail, namely the leftmost starting point, $(-r-1,1)$,
of the lattice paths, where the leftmost blue path starts. If this
endpoint could actually be reached, then the resulting object would
clearly not be of type $({\bold\TP}^g, {\bold\TP}^b)$. So we have
to show that this is impossible. (Goulden left out this
indispensable step in \cite[Theorem~2.2]{Goulden:Schur}, but we
shall close this small gap immediately.)

\begin{obs}
\label{obs:colour-changing}
The following properties of changing trails are immediate:
\begin{itemize}
\item If some edge of a changing trail is used by paths of {\em both\/}
colours green and blue, then it is necessarily traversed in both
orientations and thus forms a {\em changing trail\/} (which is ``cycle--like'')
by itself.
\item Two changing trails may well {\em touch\/}
each other (i.e., have some vertex in common), but can {\em never cross}.
\end{itemize}
\end{obs}

\begin{figure}
\caption{Illustration of the second changing trails for cases A and B.}
\label{fig:2paths}
 \begin{picture}(30,10)(0,0)
  \put(1,1){
   \begingrid{-5}{0}{5}{6}
    \gridcaption{Second changing trail, case A:}
    \makeaxes
    \makexscale
    \makeyscale
    % The recoulouring path 2:
    \thickcircle{-4}{1}
         \put(-4,1.2){\line(0,1){3.8}}
         \put(-4,5){\line(1,0){1.8}}
         \put(-2.3,4.7){\oval(0.6,0.6)[rt]}
         \put(-2,4.7){\line(0,-1){0.7}}
         \put(-2,4){\line(1,0){1.7}}
         \put(-0.3,3.7){\oval(0.6,0.6)[rt]}
         \put(0,3.7){\line(0,-1){0.7}}
         \put(0,3){\line(-1,0){1}}
         \put(-1,3){\line(0,-1){1}}
         \put(-1,2){\line(1,0){2}}
         \put(1,2){\line(0,1){2}}
         \put(1,4){\line(-1,0){0.7}}
         \put(0.3,4.3){\oval(0.6,0.6)[lb]}
         \put(0,4.3){\line(0,1){0.7}}
         \put(0,5){\vector(-1,0){1.7}}
    \thickcircle{-2}{5}
   \endgrid
  }
  \put(14,1){
   \begingrid{-5}{0}{5}{6}
    \gridcaption{Second changing trail, case B:}
    \makeaxes
    \makexscale
    \makeyscale
    % The recoulouring path 2:
    \thickcircle{-4}{1}
         \put(-4,1.2){\line(0,1){2.8}}
         \put(-4,4){\line(1,0){1.0}}
         % \put(-3.3,3.7){\oval(0.6,0.6)[rt]}
         \put(-3,4){\line(0,-1){2.0}}
         \put(-3,2){\line(1,0){1.0}}
         % \put(-2.3,1.7){\oval(0.6,0.6)[rt]}
         \put(-2,2){\line(0,-1){1.0}}
         \put(-2,1){\vector(1,0){0.8}}
    \thickcircle{-1}{1}
   \endgrid
  }
 \end{picture}
\end{figure}
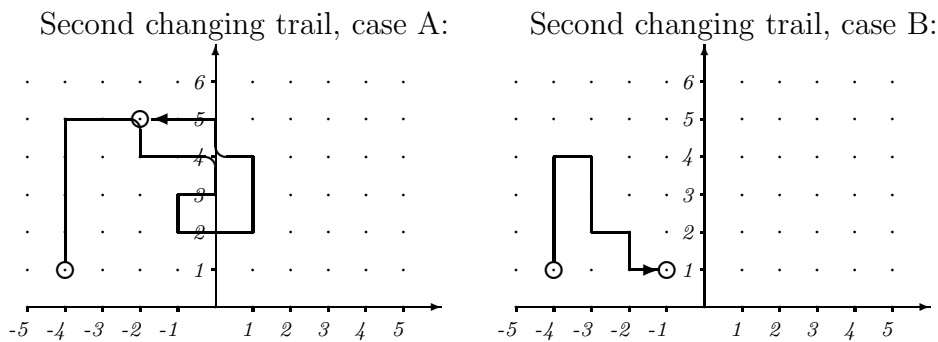

Now observe that in Case A, there is also a second possible
starting point of a ``path--like'' changing trail, namely the
left--most endpoint $(\lambda_{r+1}-r-1,N)$ of the lattice paths
(see the left picture in Figure~\ref{fig:2paths}). Likewise, in
Case B, there is a second possible starting point of a
``path--like'' changing trail, namely the rightmost starting point
$(-1,1)$ of the lattice paths (see the right picture in
Figure~\ref{fig:2paths}).

In both cases, if the changing trail starting in $(\lambda_1-1,N)$
would reach the
leftmost starting point $(-r-1,1)$ of the lattice paths, it clearly would
{\em cross\/}
this other ``path--like'' changing trail; a contradiction to
Observation~\ref{obs:colour-changing}. (The pictures in Figure~\ref{fig:2paths}
shows these other changing trails for the examples in Figures
\ref{fig:caseA} and \ref{fig:caseB}, respectively.)

This finishes the proof.
\end{proof}

%%%@@@ end input from parts/bijection

\section{The bijective construction, generalized}
\label{general}

%%%@@@ input from parts/general
It is immediately obvious that the bijective construction used in the
proof of Theorem~\ref{thm:general} is not at all restricted to the
special situation of Theorem~\ref{thm:general}: We can {\em always\/}
consider the product of two (arbitrary) skew Schur functions as
generating functions of certain ``two--coloured graphs'' derived
from the lattice path interpretation, as above. Determining the
changing trails which start in some fixed set of
starting points and recolouring their edges will {\em always\/}
yield an injective (and, clearly, weight--preserving) mapping: The
only issue which needs extra care is surjectivity.

\begin{figure}
\caption{Illustration of a {\em perfect\/} noncrossing matching in $K_{8}$.}
\label{fig:matchings}
{
\unitlength = 2.5mm
\begin{picture}(25,12)(-12,-6)
\put(4.619397662556433,
   1.913417161825448){\circle*{0.4}}
\put(1.913417161825449,
   4.619397662556433){\circle*{0.4}}
% Edge:
\put(1.913417161825449,
   4.619397662556433){\line(1,-1){2.7}}
\put(-1.913417161825448,
   4.619397662556433){\circle*{0.4}}
% Edge:
\put(-1.913417161825449,
   4.619397662556433){\line(-1,-1){2.7}}
\put(-4.619397662556433,
   1.913417161825449){\circle*{0.4}}
\put(-4.619397662556433,
   -1.913417161825449){\circle*{0.4}}
% Edge:
\put(-4.619397662556433,
   -1.913417161825449){\line(1,0){9.4}}
\put(-1.913417161825448,
   -4.619397662556433){\circle*{0.4}}
% Edge:
\put(-1.913417161825448,
   -4.619397662556433){\line(1,0){4.0}}
\put(1.913417161825449,
   -4.619397662556433){\circle*{0.4}}
\put(4.619397662556433,
   -1.913417161825448){\circle*{0.4}}\end{picture}
}
\end{figure}
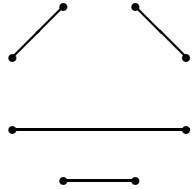

In the proof of Theorem~\ref{thm:general} we saw that the argument showing
surjectivity boils down to a very simple graph--theoretic reasoning. We
shall recast this simple reasoning into a general statement:

\begin{obs}
\label{obs:Kn}
Consider the complete graph $K_{2n}$ with $2n$ vertices,
numbered $1,2,\dots,2n$, and represent its vertices
as points on the unit circle (i.e., vertex number $m$ is represented as
$e^{2m\pi \sqrt{-1} }$); represent the edges as straight lines connecting
the corresponding vertices.
Call a matching in this graph
{\em noncrossing\/} if no two of its edges cross each other in this
geometric representation (see Figure~\ref{fig:matchings} for an illustration).
Then we have:

Any edge which belongs to a {\em perfect} noncrossing matching must
connect an odd--numbered vertex to an even--numbered vertex.
\end{obs}

\begin{rem}
Note that the number of {\em perfect\/} noncrossing matchings in $K_{2n}$
is the Catalan number $C_n$ (see \cite[p.~222]{Stanley2}).
% according to Stanley \cite[p.~222]{Stanley2},
% this result is apparently due to A.~Errera \cite{Errera}.
\end{rem}

\begin{rem}
Note that the argument proving surjectivity in Theorem~\ref{thm:general}
amounts to the fact that the two possible ``path--like'' changing
trails connecting the four possible starting points and end points
$(\lambda_1-1,N)$, $(\lambda_{r+1}-r-1,N)$, $(-r-1,1)$ and $(-1,1)$ must
correspond to a {\em noncrossing\/} perfect matching of the complete graph
$K_4$. 
\end{rem}

We shall derive a general statement for skew Schur functions:

Let $\lambda=(\lambda_1,\dots,\lambda_r)$ and
$\mu=(\mu_1,\dots,\mu_r)$ be partitions with
$0\leq\mu_i\leq\lambda_i$ for $1\leq i\leq r$; let
$\sigma=(\sigma_1,\dots,\sigma_r)$ and $\tau=(\tau_1,\dots,\tau_r)$
be partitions with $0\leq\tau_i\leq\sigma_i$ for
$1\leq i\leq r$.

\begin{rem}
\label{rem:zero-parts}
We intentionally allow parts of length 0 in the partitions $\lambda$
and $\sigma$: This is equivalent to allowing them to have
different numbers of parts.
\end{rem}

Interpret $s_{\lambda/\mu}$ as the generating function of the
family of nonintersecting
lattice paths $(P^b_1,\dots,P^b_r)$, where
$P^b_i$ starts at $(\mu_{i}-i, 1)$ and ends at $(\lambda_{i}-i, N)$.
Colour the corresponding lattice paths blue.

Interpret $s_{\sigma/\tau}$ as the generating function of the family of nonintersecting
lattice paths $(P^g_1,\dots,P^g_r)$, where $P^g_i$ starts at
$(\tau_{i}+t-i, 1)$
and ends at $(\sigma_{i}+t-i, N)$. Colour the corresponding lattice paths
green. Here, $t$ is an arbitrary but fixed integer which indicates the
horizontal offset of the green paths with respect to the blue paths.

Consider the sequence of possible starting points of ``path--like''
changing trails of the corresponding two-coloured graph, in the
sense of Section~\ref{bijection}, where the end--points of the
lattice paths appear in order from right to left in this sequence,
followed by the starting points of the lattice paths in order from
left to right. Note that the number of such points is even, $2k$,
say. More precisely, consider $(x_1,N),\dots (x_l,N),$ followed by
$(x_{l+1},1),\dots,(x_{2k},1)$, where
\begin{multline*}
\{x_1,\dots,x_l\}=
\{i:\;\lambda_{i}-i\neq\sigma_{j}+t-j\text{ for }1\leq j\leq r\} \cup\\
\{j:\;\sigma_{j}+t-j\neq\lambda_{i}-i\text{ for }1\leq i\leq r\},
\end{multline*}
$x_1>x_2>\dots >x_l$, and where
\begin{multline*}
\{x_{l+1},\dots,x_{2k}\}=
\{i:\;\mu_{i}-i\neq\tau_{j}+t-j\text{ for }1\leq j\leq r\} \cup\\
\{j:\;\tau_{j}+t-j\neq\mu_{i}-i\text{ for }1\leq i\leq r\},
\end{multline*}
$x_{l+1}<\dots <x_{2k}.$

Denote this sequence of points $(x_i,.)$ by $(Q_i),1\leq i\leq 2k$.
For $1\leq i\leq l$, blue points $Q_i$ are coloured black and green
points $Q_i$ are coloured white. For $l+1\leq i\leq 2k$, blue
points $Q_i$ are coloured white and green points $Q_i$ are coloured
black. Points with even index are called even, points with odd
index are called odd. Then the following lemma is immediate:

\begin{lem}
\label{lem:parity}
A path--like changing trail in the two--coloured graph defined above can only
connect points of different colours (out of black and white) {\em
and\/} of different parity (by Observation \ref{obs:Kn});
e.g., some white $Q_{2m}$ and some
black $Q_{2n+1}$.
\end{lem}

Now fix an arbitrary subset of points $\{Q_{i_1},\dots,Q_{i_m}\}$.
Start with an arbitrary two--coloured graph from
$s_{\lambda/\mu}s_{\sigma/\tau}$ (interpreted again as the product of the
generating functions of the corresponding families of
nonintersecting lattice paths) and recolour the changing trails
starting in $Q_{i_1},\dots,Q_{i_m}$. In general, this will give
another two--coloured graph, which can be interpreted as belonging
to some other
$s_{\lambda^\prime/\mu^\prime}s_{\sigma^\prime/\tau^\prime}$. Take
an arbitrary object (i.e., two--coloured graph) from
$s_{\lambda^\prime/\mu^\prime}s_{\sigma^\prime/\tau^\prime}$ and
repeat the same recolouring operation as long as it generates some
``new'' (yet unseen) object.

The set of objects thus generated decomposes into two disjoint
sets: One set, $O_0$, encompasses all objects which show the same
colouring of points $Q_{i_1},\dots,Q_{i_m}$ as in the starting
object; the other, $O_1$ encompasses the objects with the opposite
colouring for these points.

It is clear that recolouring changing trails which start in points
$Q_{i_1},\dots,Q_{i_m}$ establishes a bijection between $O_0$ and
$O_1$.

On the other hand, each object in $O_0$ belongs to some
$s_{\lambda^\pprime/\mu^\pprime}s_{\sigma^\pprime/\tau^\pprime}$:
Denote the set of all the corresponding quadruples
$(\lambda^\pprime,\mu^\pprime,\sigma^\pprime,\tau^\pprime)$ which
occur in this sense by $S_0$. The same consideration applies to
$O_1$: Denote by $S_1$ the corresponding set  of
quadruples $(\lambda^\prime,\mu^\prime,\sigma^\prime,\tau^\prime)$.

\begin{lem}
\label{lem:most-general}
Given the above definitions, we have the following ``generic''
identity for skew Schur functions:
\begin{equation}
\label{eq:most-general}
\sum_{(\lambda^\prime,\mu^\prime,\sigma^\prime,\tau^\prime)\in S_1}
                s_{\lambda^\prime/\mu^\prime}s_{\sigma^\prime/\tau^\prime}=
\sum_{(\lambda^\pprime,\mu^\pprime,\sigma^\pprime,\tau^\pprime)\in S_0}
                s_{\lambda^\pprime/\mu^\pprime}s_{\sigma^\pprime/\tau^\pprime}.
\end{equation}
\end{lem}

This statement is certainly as general as useless: Let us
specialize to a somewhat ``friendlier'' assertion.
\begin{lem}
\label{lem:general-skew}
Given the above definitions, assume that all black points have the same
parity, and that all white points have the same parity. Then \eqref{eq:most-general} specializes to
\begin{equation}
s_{\lambda/\mu}s_{\sigma/\tau}=
\sum_{(\lambda^\prime,\mu^\prime,\sigma^\prime,\tau^\prime)\in S_1}
                s_{\lambda^\prime/\mu^\prime}s_{\sigma^\prime/\tau^\prime},
\end{equation}
where $S_1$ encompasses all the quadruples
$(\lambda^\prime,\mu^\prime,\sigma^\prime,\tau^\prime)$
which correspond to any two--coloured graph object that can be obtained
by recolouring the changing trails starting in points
$Q_{i_1},\dots,Q_{i_m}$ in any ``initial'' two--coloured graph object
from $s_{\lambda/\mu}s_{\sigma/\tau}$.
\end{lem}
\begin{proof}
Without loss of generality we may assume that all even points are
white and all odd points are black in
$s_{\lambda/\mu}s_{\sigma/\tau}$. By recolouring changing trails,
all the points $Q_{i_1},\dots,Q_{i_m}$ are matched with points of
opposite colour and parity.

So if $Q_i$ is odd and black, then the recolouring trail starting at
$Q_i$ connects it which some other point $Q_k$ which is even and white:
After recolouring, $Q_i$ is odd and white, and the recolouring operation
altogether yields some two--coloured graph object belonging to some
$s_{\lambda^\prime/\mu^\prime}s_{\sigma^\prime/\tau^\prime}$.

Now if we apply the recolouring operation to an {\em arbitrary\/}
object from
$s_{\lambda^\prime/\mu^\prime}s_{\sigma^\prime/\tau^\prime}$, the
only possible partners for ``wrongly--coloured'' $Q_i$ (odd, but
white) is another ``wrongly--coloured'' $Q_j$ (even, but black).
Hence this operation takes objects from
$s_{\lambda^\prime/\mu^\prime}s_{\sigma^\prime/\tau^\prime}$ back
to $s_{\lambda/\mu}s_{\sigma/\tau}$.
\end{proof}

%%%@@@ end input from parts/general

\section{Proof of the Pl\"ucker relations}
\label{pluecker}

%%%@@@ input from parts/pluecker
In order to show that the general assertions of
Section~\ref{general} do in fact lead to some interesting
identities, we give a proof of the Pl\"ucker relations
(Theorem~\ref{thm:pluecker}), which is based on
Lemma~\ref{lem:general-skew}.

\begin{proof}
In the notation of Section~\ref{general}, let
$$\lambda=\left(2n(n-1),2(n-1)^2,\dots,4(n-1),2(n-1)\right)$$
and
$$\sigma=\left((2n-1)(n-1),(2n-3)(n-1),\dots, 3(n-1),n-1\right),$$
$\mu=\tau=(0,\dots,0)$; and choose horizontal offset $t=0$. I.e., interpret
$s_\lambda s_\sigma$ as the generating function of two--coloured
graph objects consisting of two $n$-tuples of nonintersecting
lattice paths, coloured
green and blue, respectively,
where green path $P^g_i$ starts at $(-i,1)$ and ends at
$(\lambda_i-i,N)$, and where blue path $P^b_i$ starts at $(-i,1)$
and ends at $(\sigma_i-i,N)$.

Observe that this setting obeys the assumptions of
Lemma~\ref{lem:general-skew}.

Now consider the set of green endpoints $\{Q_1,\dots,Q_k\}$, where
$Q_i=(\lambda_{r_i}-r_i,N)$.
(Here, $1\leq r_1< r_2<\dots< r_k\leq n$ is the fixed
list of integers from Theorem~\ref{thm:pluecker}.)
Recolouring changing trails which
start at these points amounts to determining the set
$\{R_1,\dots,R_k\}$ of respective
endpoints of the changing trails, and changing colours.

Assume that $R_i=(\sigma_{s_i}-s_i,N)$, then in terms of the associated
Schur functions Lemma~\ref{lem:general-skew} directly leads to the identity:
\begin{equation}
\label{eq:schur-pluecker}
s_\lambda s_\sigma = \sum_{1\leq s_1<\dots< s_k\leq n}
        s_{(\lambda_1,\dots,\sigma_{s_1},\dots,\sigma_{s_k},\dots,\lambda_n)}
        s_{(\sigma_1,\dots,\lambda_{r_1},\dots,\lambda_{r_k},\dots,\sigma_n)},
\end{equation}
where the notation of the summands means that parts $\lambda_{r_i}$
were exchanged with parts $\sigma_{s_i}$, respectively.

By the Jacobi--Trudi identity \eqref{eq:Jacobi-Trudi} and
Observation~\ref{obs:christian},
\eqref{eq:pluecker} and \eqref{eq:schur-pluecker} are in fact equivalent.
\end{proof}

\begin{rem}
In fact, even the quite general assertion of Lemma~\ref{lem:most-general}
can be generalized further: So far, our lattice paths always had
starting points and end points at the same horizontal lines $(.,1)$ and
$(.,N)$, corresponding to the range of variables $x_1,\dots,x_N$. Dropping
this constraint yields Schur functions with different ranges of variables
(e.g., $s_\lambda(x_4,x_5,x_6)$). Recalling that (see
Remark~\ref{rem:zero-parts}) we actually also do allow partitions of different
lengths, it is easy to see that Theorem~5 in \cite{Krattenthaler}
(which is a generalization of Ciucu's Schur function identity
\eqref{eq:ciucu}) can be proved in the same way as
Lemma~\ref{lem:most-general}.
\end{rem}

%%%@@@ end input from parts/pluecker

% markus
\section{Kleber's Theorem}
\label{kleber}

%%%@@@ input from parts/kleber
% markus
The theorem \cite[Thm.~3.2]{Kleber} is expressed in terms of
certain operations on Ferrers boards (called Young
diagrams in \cite{Kleber}): In order to state it, we need to describe the
relevant notation.

First, we introduce a particular way of drawing the Ferrers board
of $\lambda=(\lambda_1,\dots,\lambda_r)$ in the plane: Let
$x_1>x_2>\dots>x_n>x_{n+1}=0$ be the ordered list of the {\em
distinct\/} parts contained in $\lambda$; set $y_i=\text{the number
of parts of }\lambda\text{ which are }\geq x_i$.

Setting $y_0=0$, we have $0=y_0<y_1<\dots<y_n$, and $(x_i),(y_i)$
simply yield another encoding of the partition $\lambda$:
$$\lambda=(x_1^{y_1-y_0},x_2^{y_2-y_1},\dots,x_n^{y_n-y_{n-1}}).$$

Now consider the $n$ points $(x_1,-y_1),(x_2,-y_2),\dots,(x_n,-y_n)$ in
the plane: The Ferrers board of $\lambda$ is represented as the set of
points $(x,-y)$ such that:
\begin{align*}
x\geq 0\text{ and } y\geq 0,&\\
x\leq x_i\text{ and }y\leq y_i&\text{ for some } i.
\end{align*}

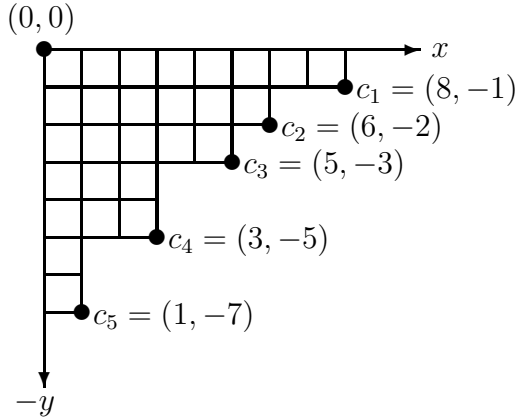
\begin{figure}
\caption{Illustration of outer corners and special drawing of Ferrers
        board for partition $\lambda=(8,6,5,3,3,1,1)$.}
\label{fig:outer-corners}
 \begin{picture}(25,12)(0,0)
  % \graphpaper[5](0,0)(25,12)
  \put(8,0){
   \begin{picture}(10,10)(0,-10)
         % Horizontal lines
         \put(0,0){\vector(1,0){10}}    \put(10.3,-0.2){\hbox{$x$}}
         \put(0,-1){\line(1,0){8}}
         \put(0,-2){\line(1,0){6}}
         \put(0,-3){\line(1,0){5}}
         \put(0,-4){\line(1,0){3}}
         \put(0,-5){\line(1,0){3}}
         \put(0,-6){\line(1,0){1}}
         \put(0,-7){\line(1,0){1}}
         % Vertical lines
         \put(0,0){\vector(0,-1){9}}    \put(-0.8,-9.6){\hbox{$-y$}}
         \put(1,0){\line(0,-1){7}}
         \put(2,0){\line(0,-1){5}}
         \put(3,0){\line(0,-1){5}}
         \put(4,0){\line(0,-1){3}}
         \put(5,0){\line(0,-1){3}}
         \put(6,0){\line(0,-1){2}}
         \put(7,0){\line(0,-1){1}}
         \put(8,0){\line(0,-1){1}}
         % Outer corners
         \put(8,-1){\circle*{0.4}}              \put(8.3,-1.3){\hbox{$c_1=(8,-1)$}}
         \put(6,-2){\circle*{0.4}}              \put(6.3,-2.3){\hbox{$c_2=(6,-2)$}}
         \put(5,-3){\circle*{0.4}}              \put(5.3,-3.3){\hbox{$c_3=(5,-3)$}}
         \put(3,-5){\circle*{0.4}}              \put(3.3,-5.3){\hbox{$c_4=(3,-5)$}}
         \put(1,-7){\circle*{0.4}}              \put(1.3,-7.3){\hbox{$c_5=(1,-7)$}}
         % Origin
         \put(0,0){\circle*{0.4}}               \put(-1,0.6){\hbox{$(0,0)$}}
        \end{picture}
  }
 \end{picture}
\end{figure}

Figure \ref{fig:outer-corners} illustrates this concept. The $n$ points
$c_1=(x_1,-y_1),c_2=(x_2,-y_2),\dots,c_n=(x_n,-y_n)$ are called {\em outside corners\/},
the $n+1$ points $(x_1,-y_0),(x_2,-y_1),\dots,(x_{n+1},-y_n)$
are called {\em inside corners\/}.

Now we are in a position to define two operations on partitions:
In the above notation, take two integers $i,j$ such that $1\leq i\leq j\leq n$
and define two partitions derived from the original $\lambda$
via manipulating the inside and outside corners of its
associated Ferrers board:
\begin{align*}
\pi^i_j(\lambda):&\text{ add $1$ to each of }x_{i+1},\dots,x_j;y_i,\dots,y_j,\\
\mu^i_j(\lambda):&\text{ add $-1$ to each of }x_{i+1},\dots,x_j;y_i,\dots,y_j.
\end{align*}

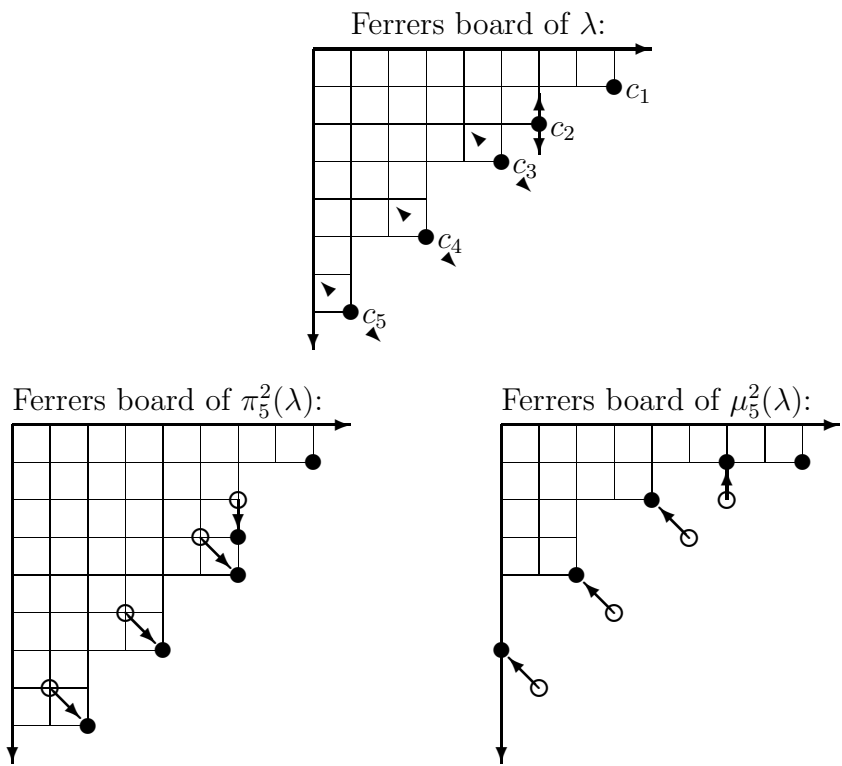
\begin{figure}
\caption{Illustration of operations $\pi_j^i$ and $\mu_j^i$ for
        $i=2$, $j=5$ applied to $\lambda=(8,6,5,3,3,1,1)$.}
\label{fig:kleber}
 \begin{picture}(25,21)(0,-1)
  % \graphpaper[5](0,0)(25,19)
  \linethickness{0.01pt}
  \put(10,10){
   \begin{picture}(10,8)(0,-8)
         % Heading
         \put(1,0.4){\hbox{Ferrers board of $\lambda$:}}
         % Axes
         \thicklines
         \put(0,0){\vector(1,0){9}}
         \put(0,0){\vector(0,-1){8}}
    \linethickness{0.01pt}
         % Horizontal lines
         \put(0,0){\line(1,0){8}}
         \put(0,-1){\line(1,0){8}}
         \put(0,-2){\line(1,0){6}}
         \put(0,-3){\line(1,0){5}}
         \put(0,-4){\line(1,0){3}}
         \put(0,-5){\line(1,0){3}}
         \put(0,-6){\line(1,0){1}}
         \put(0,-7){\line(1,0){1}}
         % Vertical lines
         \put(0,0){\line(0,-1){7}}
         \put(1,0){\line(0,-1){7}}
         \put(2,0){\line(0,-1){5}}
         \put(3,0){\line(0,-1){5}}
         \put(4,0){\line(0,-1){3}}
         \put(5,0){\line(0,-1){3}}
         \put(6,0){\line(0,-1){2}}
         \put(7,0){\line(0,-1){1}}
         \put(8,0){\line(0,-1){1}}
         % Outer corners
         \thicklines
         \put(8,-1){\circle*{0.4}}              \put(8.3,-1.3){\hbox{$c_1$}}
         \put(6,-2){\circle*{0.4}}              \put(6.3,-2.3){\hbox{$c_2$}}
         \put(6,-2.2){\vector(0,-1){0.6}}
         \put(6,-1.8){\vector(0,1){0.6}}
         \put(5,-3){\circle*{0.4}}              \put(5.3,-3.3){\hbox{$c_3$}}
         \put(5.2,-3.2){\vector(1,-1){0.6}}
         \put(4.8,-2.8){\vector(-1,1){0.6}}
         \put(3,-5){\circle*{0.4}}              \put(3.3,-5.3){\hbox{$c_4$}}
         \put(3.2,-5.2){\vector(1,-1){0.6}}
         \put(2.8,-4.8){\vector(-1,1){0.6}}
         \put(1,-7){\circle*{0.4}}              \put(1.3,-7.3){\hbox{$c_5$}}
         \put(1.2,-7.2){\vector(1,-1){0.6}}
         \put(0.8,-6.8){\vector(-1,1){0.6}}
        \end{picture}
  }
  \put(2,0){
   \begin{picture}(10,8)(0,-8)
         % Heading
         \put(0,0.4){\hbox{Ferrers board of $\pi_5^2(\lambda)$:}}
         % Axes
         \thicklines
         \put(0,0){\vector(1,0){9}}
         \put(0,0){\vector(0,-1){9}}
    \linethickness{0.01pt}
         % Horizontal lines
         \put(0,0){\line(1,0){8}}
         \put(0,-1){\line(1,0){8}}
         \put(0,-2){\line(1,0){6}}
         \put(0,-3){\line(1,0){6}}
         \put(0,-4){\line(1,0){6}}
         \put(0,-5){\line(1,0){4}}
         \put(0,-6){\line(1,0){4}}
         \put(0,-7){\line(1,0){2}}
         \put(0,-8){\line(1,0){2}}
         % Vertical lines
         \put(0,0){\line(0,-1){8}}
         \put(1,0){\line(0,-1){8}}
         \put(2,0){\line(0,-1){8}}
         \put(3,0){\line(0,-1){6}}
         \put(4,0){\line(0,-1){6}}
         \put(5,0){\line(0,-1){4}}
         \put(6,0){\line(0,-1){4}}
         \put(7,0){\line(0,-1){1}}
         \put(8,0){\line(0,-1){1}}
         % Outer corners
         \thicklines
         \put(8,-1){\circle*{0.4}}
         \put(6,-2){\circle{0.4}}
         \put(6,-2){\vector(0,-1){0.8}}
         \put(6,-3){\circle*{0.4}}
         \put(5,-3){\circle{0.4}}
         \put(5,-3){\vector(1,-1){0.8}}
         \put(6,-4){\circle*{0.4}}
         \put(3,-5){\circle{0.4}}
         \put(3,-5){\vector(1,-1){0.8}}
         \put(4,-6){\circle*{0.4}}
         \put(1,-7){\circle{0.4}}
         \put(1,-7){\vector(1,-1){0.8}}
         \put(2,-8){\circle*{0.4}}
        \end{picture}
  }
  \put(15,0){
   \begin{picture}(10,8)(0,-8)
         % Heading
         \put(0,0.4){\hbox{Ferrers board of $\mu_5^2(\lambda)$:}}
         % Axes
         \thicklines
         \put(0,0){\vector(1,0){9}}
         \put(0,0){\vector(0,-1){9}}
    \linethickness{0.01pt}
         % Horizontal lines
         \put(0,0){\line(1,0){8}}
         \put(0,-1){\line(1,0){8}}
         \put(0,-2){\line(1,0){4}}
         \put(0,-3){\line(1,0){2}}
         \put(0,-4){\line(1,0){2}}
         % Vertical lines
         \put(0,0){\line(0,-1){4}}
         \put(1,0){\line(0,-1){4}}
         \put(2,0){\line(0,-1){4}}
         \put(3,0){\line(0,-1){2}}
         \put(4,0){\line(0,-1){2}}
         \put(5,0){\line(0,-1){1}}
         \put(6,0){\line(0,-1){1}}
         \put(7,0){\line(0,-1){1}}
         \put(8,0){\line(0,-1){1}}
         % Outer corners
         \thicklines
         \put(8,-1){\circle*{0.4}}%             \put(8.3,-1.3){\hbox{$c_1$}}
         \put(6,-2){\circle{0.4}}%              \put(6.3,-2.3){\hbox{$c_2$}}
         \put(6,-2){\vector(0,1){0.8}}
         \put(6,-1){\circle*{0.4}}%             \put(6.3,-2.3){\hbox{$c_2$}}
         \put(5,-3){\circle{0.4}}%              \put(6.3,-2.3){\hbox{$c_2$}}
         \put(5,-3){\vector(-1,1){0.8}}
         \put(4,-2){\circle*{0.4}}%             \put(5.3,-3.3){\hbox{$c_3$}}
         \put(3,-5){\circle{0.4}}%              \put(6.3,-2.3){\hbox{$c_2$}}
         \put(3,-5){\vector(-1,1){0.8}}
         \put(2,-4){\circle*{0.4}}%             \put(3.3,-5.3){\hbox{$c_4$}}
         \put(1,-7){\circle{0.4}}
         \put(1,-7){\vector(-1,1){0.8}}
         \put(0,-6){\circle*{0.4}}
        \end{picture}
  }
 \end{picture}
\end{figure}

These operations add or remove, respectively, a {\em border strip\/}
that reaches from the $i$-th outside corner to the $j$-th inside corner
(see Figure \ref{fig:kleber}).

We need to add or remove {\em nested border strips\/}: Given integers
$1\leq i_1<\dots<i_k\leq j_k<\dots j_1\leq n$, we define
\begin{align*}
\pi^{i_1,\dots,i_k}_{j_1,\dots,j_k} & =
        \pi^{i_1}_{j_1}\circ \dots\circ\pi^{i_k}_{j_k},\\
\mu^{i_1,\dots,i_k}_{j_1,\dots,j_k} & =
        \mu^{i_1}_{j_1}\circ \dots\circ\mu^{i_k}_{j_k}.
\end{align*}

Note that the corners which are shifted by these operations might
not appear as corners in the geometric sense any more; nevertheless
we consider them as the object for subsequent operations $\pi$ and $\mu$:
Nesting $\pi$ and $\mu$ in this sense yields something which can be interpreted
again as a partition, since we always have $x_i\geq x_{i+1}$
and $y_i\leq y_{i+1}$ (see Figure~\ref{fig:kleber}.)

\begin{figure}
\caption{Illustration of operation $\lambda\pm\omega_l$ for
        $l=y_4=5$, applied to $\lambda=(8,6,5,3,3,1,1)$.}
\label{fig:kleber2}
 \begin{picture}(25,11)(0,-1)
  % \graphpaper[5](0,0)(25,19)
  \linethickness{0.01pt}
  \put(2,0){
   \begin{picture}(10,8)(0,-8)
         % Heading
         \put(0,0.4){\hbox{Ferrers board of $\lambda+\omega_5$:}}
         % Axes
         \thicklines
         \put(0,0){\vector(1,0){10}}
         \put(0,0){\vector(0,-1){9}}
    \linethickness{0.01pt}
         % Horizontal lines
         \put(0,0){\line(1,0){9}}
         \put(0,-1){\line(1,0){9}}
         \put(0,-2){\line(1,0){7}}
         \put(0,-3){\line(1,0){6}}
         \put(0,-4){\line(1,0){4}}
         \put(0,-5){\line(1,0){4}}
         \put(0,-6){\line(1,0){1}}
         \put(0,-7){\line(1,0){1}}
         % Vertical lines
         \put(0,0){\line(0,-1){7}}
         \put(1,0){\line(0,-1){7}}
         \put(2,0){\line(0,-1){5}}
         \put(3,0){\line(0,-1){5}}
         \put(4,0){\line(0,-1){5}}
         \put(5,0){\line(0,-1){3}}
         \put(6,0){\line(0,-1){3}}
         \put(7,0){\line(0,-1){2}}
         \put(8,0){\line(0,-1){1}}
         \put(9,0){\line(0,-1){1}}
         % Outer corners
         \thicklines
         \put(8,-1){\vector(1,0){0.8}}
         \put(9,-1){\circle*{0.4}}
         \put(8,-1){\circle{0.4}}
         \put(6,-2){\vector(1,0){0.8}}
         \put(7,-2){\circle*{0.4}}
         \put(6,-2){\circle{0.4}}
         \put(5,-3){\vector(1,0){0.8}}
         \put(6,-3){\circle*{0.4}}
         \put(5,-3){\circle{0.4}}
         \put(3,-5){\vector(1,0){0.8}}
         \put(4,-5){\circle*{0.4}}
         \put(3,-5){\circle{0.4}}
         \put(1,-7){\circle*{0.4}}
        \end{picture}
  }
  \put(15,0){
   \begin{picture}(10,8)(0,-8)
         % Heading
         \put(0,0.4){\hbox{Ferrers board of $\lambda-\omega_5$:}}
         % Axes
         \thicklines
         \put(0,0){\vector(1,0){9}}
         \put(0,0){\vector(0,-1){9}}
    \linethickness{0.01pt}
         % Horizontal lines
         \put(0,0){\line(1,0){7}}
         \put(0,-1){\line(1,0){7}}
         \put(0,-2){\line(1,0){5}}
         \put(0,-3){\line(1,0){4}}
         \put(0,-4){\line(1,0){2}}
         \put(0,-5){\line(1,0){2}}
         \put(0,-6){\line(1,0){1}}
         \put(0,-7){\line(1,0){1}}
         % Vertical lines
         \put(0,0){\line(0,-1){7}}
         \put(1,0){\line(0,-1){7}}
         \put(2,0){\line(0,-1){5}}
         \put(3,0){\line(0,-1){3}}
         \put(4,0){\line(0,-1){3}}
         \put(5,0){\line(0,-1){2}}
         \put(6,0){\line(0,-1){1}}
         \put(7,0){\line(0,-1){1}}
         % Outer corners
         \thicklines
         \put(8,-1){\vector(-1,0){0.8}}
         \put(7,-1){\circle*{0.4}}
         \put(8,-1){\circle{0.4}}
         \put(6,-2){\vector(-1,0){0.8}}
         \put(5,-2){\circle*{0.4}}
         \put(6,-2){\circle{0.4}}
         \put(5,-3){\vector(-1,0){0.8}}
         \put(4,-3){\circle*{0.4}}
         \put(5,-3){\circle{0.4}}
         \put(3,-5){\vector(-1,0){0.8}}
         \put(2,-5){\circle*{0.4}}
         \put(3,-5){\circle{0.4}}
         \put(1,-7){\circle*{0.4}}
        \end{picture}
  }
 \end{picture}
\end{figure}
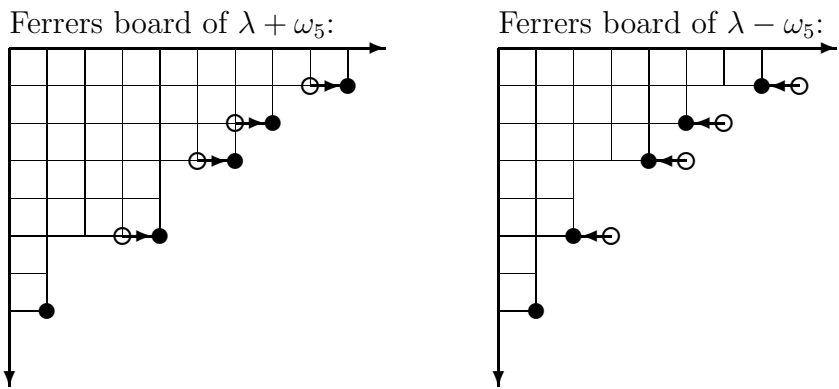

The last operation we need is the following: In the above notation,
let $k$ be an integer, $1\leq k\leq n$. Clearly, the Ferrers board
contains at least one column of length $l=y_k$: Adding or removing
some column of length $l$ amounts to adding $\pm 1$ to all
coordinates $x_i$, $1\leq i\leq k$. We denote this operation by
$\lambda\pm\omega_l$. (See Figure~\ref{fig:kleber2}.)

\begin{thm}[Theorem 3.2 in \cite{Kleber}]
\label{thm:Kleber}
Let $\lambda=(\lambda_1,\lambda_2,\dots,\lambda_r)$ be a partition
with $n$ outside corners. For an arbitrary integer $k$, $1\leq
k\leq n$, set $l=y_k$ (in the above notation). Then we have:
\begin{multline}
\label{eq:kleber}
s_\lambda s_\lambda = \\
s_{\lambda+\omega_l}s_{\lambda-\omega_l} +
\sum_{m\geq 1}
        \sum_{\substack{
                        1\leq i_1<\dots<i_m\leq k\\
                        k\leq j_m<\dots<j_1\leq n
                }}
                (-1)^{m-1}s_{\pi^{i_1,\dots,i_m}_{j_1,\dots,j_m}(\lambda)}
                          s_{\mu^{i_1,\dots,i_m}_{j_1,\dots,j_m}(\lambda)}.
\end{multline}
\end{thm}

The connections between Ferrers boards and nonintersecting lattice
paths were illustrated in Section~\ref{background}: Here we have to
give the proper ``translation'' of operations
$\pi_{j}^{i}$ and $\mu_{j}^{i}$ to nonintersecting lattice paths.

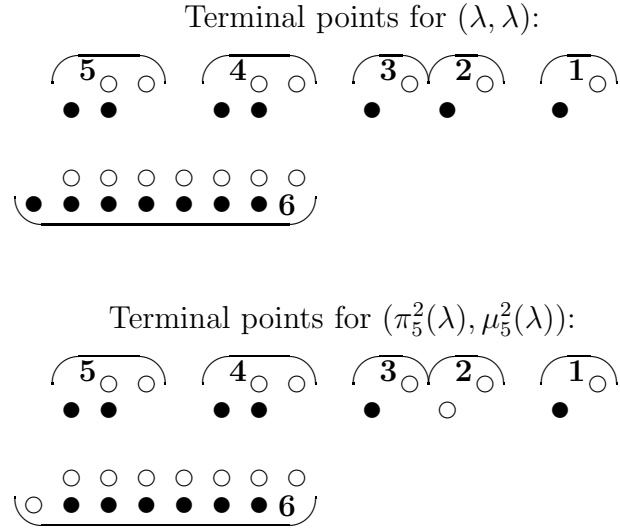
\begin{figure}
\caption{Illustration of operations $\pi_j^i$ and $\mu_j^i$ for
        $i=2$, $j=4$ applied to $\lambda=(8,6,5,3,3,1,1)$,
        translated to lattice paths.}
\label{fig:kleber-lattice}
 \begin{picture}(25,16)(0,0)
  \linethickness{0.01pt}
  % \graphpaper[5](0,0)(25,19)
  \thinlines
  \put(4,8){
   \begin{picture}(15,6)(0,-1)
         % Heading
         \put(4,5){\hbox{Terminal points for $(\lambda,\lambda)$:}}
         % Starting points, green
         \put(1,1){\circle{0.4}}
         \put(2,1){\circle{0.4}}
         \put(3,1){\circle{0.4}}
         \put(4,1){\circle{0.4}}
         \put(5,1){\circle{0.4}}
         \put(6,1){\circle{0.4}}
         \put(7,1){\circle{0.4}}
         % Starting points, blue
         \put(0,0.3){\circle*{0.4}}
         \put(1,0.3){\circle*{0.4}}
         \put(2,0.3){\circle*{0.4}}
         \put(3,0.3){\circle*{0.4}}
         \put(4,0.3){\circle*{0.4}}
         \put(5,0.3){\circle*{0.4}}
         \put(6,0.3){\circle*{0.4}}
         % End points, green
         \put(2,3.5){\circle{0.4}}
         \put(3,3.5){\circle{0.4}}
         \put(6,3.5){\circle{0.4}}
         \put(7,3.5){\circle{0.4}}
         \put(10,3.5){\circle{0.4}}
         \put(12,3.5){\circle{0.4}}
         \put(15,3.5){\circle{0.4}}
         % End points, blue
         \put(1,2.8){\circle*{0.4}}
         \put(2,2.8){\circle*{0.4}}
         \put(5,2.8){\circle*{0.4}}
         \put(6,2.8){\circle*{0.4}}
         \put(9,2.8){\circle*{0.4}}
         \put(11,2.8){\circle*{0.4}}
         \put(14,2.8){\circle*{0.4}}
         % Blocks
         \put(3.5,0.5){\oval(8,1.5)[b]} \put(6.5,0){\hbox{\bf 6}}
         \put(2,3.5){\oval(3,1.5)[t]}           \put(1.2,3.6){\hbox{\bf 5}}
         \put(6,3.5){\oval(3,1.5)[t]}           \put(5.2,3.6){\hbox{\bf 4}}
         \put(9.5,3.5){\oval(2,1.5)[t]} \put(9.2,3.6){\hbox{\bf 3}}
         \put(11.5,3.5){\oval(2,1.5)[t]}        \put(11.2,3.6){\hbox{\bf 2}}
         \put(14.5,3.5){\oval(2,1.5)[t]}        \put(14.2,3.6){\hbox{\bf 1}}
        \end{picture}
  }
  \put(4,0){
   \begin{picture}(15,6)(0,-1)
         % Heading
         \put(2,5){\hbox{Terminal points for $(\pi_5^2(\lambda),\mu_5^2(\lambda))$:}}
         % Starting points, green
         \put(1,1){\circle{0.4}}
         \put(2,1){\circle{0.4}}
         \put(3,1){\circle{0.4}}
         \put(4,1){\circle{0.4}}
         \put(5,1){\circle{0.4}}
         \put(6,1){\circle{0.4}}
         \put(7,1){\circle{0.4}}
         % Starting points, blue
         \put(0,0.3){\circle{0.4}}
         \put(1,0.3){\circle*{0.4}}
         \put(2,0.3){\circle*{0.4}}
         \put(3,0.3){\circle*{0.4}}
         \put(4,0.3){\circle*{0.4}}
         \put(5,0.3){\circle*{0.4}}
         \put(6,0.3){\circle*{0.4}}
         % End points, green
         \put(2,3.5){\circle{0.4}}
         \put(3,3.5){\circle{0.4}}
         \put(6,3.5){\circle{0.4}}
         \put(7,3.5){\circle{0.4}}
         \put(10,3.5){\circle{0.4}}
         \put(12,3.5){\circle{0.4}}
         \put(15,3.5){\circle{0.4}}
         % End points, blue
         \put(1,2.8){\circle*{0.4}}
         \put(2,2.8){\circle*{0.4}}
         \put(5,2.8){\circle*{0.4}}
         \put(6,2.8){\circle*{0.4}}
         \put(9,2.8){\circle*{0.4}}
         \put(11,2.8){\circle{0.4}}
         \put(14,2.8){\circle*{0.4}}
         % Blocks
         \put(3.5,0.5){\oval(8,1.5)[b]} \put(6.5,0){\hbox{\bf 6}}
         \put(2,3.5){\oval(3,1.5)[t]}           \put(1.2,3.6){\hbox{\bf 5}}
         \put(6,3.5){\oval(3,1.5)[t]}           \put(5.2,3.6){\hbox{\bf 4}}
         \put(9.5,3.5){\oval(2,1.5)[t]} \put(9.2,3.6){\hbox{\bf 3}}
         \put(11.5,3.5){\oval(2,1.5)[t]}        \put(11.2,3.6){\hbox{\bf 2}}
         \put(14.5,3.5){\oval(2,1.5)[t]}        \put(14.2,3.6){\hbox{\bf 1}}
        \end{picture}
  }
 \end{picture}
\end{figure}

First observe that the outside corners of a partition correspond to {\em
blocks of consecutive endpoints\/} (here, consecutive means ``having distance
1 in the horizontal direction'') in the lattice path interpretation: Number
these blocks from right to left by $1,2,\dots,n$, and denote the additional
block of (consecutive) starting points by $n+1$ (see the upper picture
in Figure~\ref{fig:kleber-lattice}).

Interpret some object from $s_\lambda s_\lambda$ in the same way as in
Section~\ref{general}. More precisely, let $\sigma=\lambda$, $\mu=\tau=0$
and horizontal offset $t=1$ in the general definitions preceding
Lemma~\ref{lem:most-general}. Figure~\ref{fig:kleber-lattice}
illustrates the position of starting points and end points of the corresponding
lattice paths: Blue points are drawn as black dots,
green points are drawn as white dots; blocks are indicated by horizontal
braces.

It is easy to see that the simultaneous
application of $\pi_j^i$ to the ``green object'' and of $\mu_j^i$ to the
``blue
object'' amounts to interchanging colours of the leftmost point in blue block
$i$ and of the appropriate endpoint of a corresponding
changing trail in block $j+1$ (i.e., the rightmost point in green block $j+1$
if $j<n$, or the the leftmost point in blue block $n+1$ if $j=n$; see Figure
\ref{fig:kleber-lattice}).

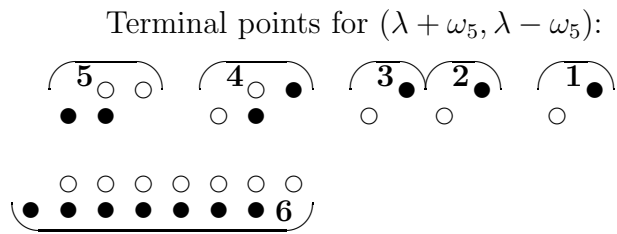
\begin{figure}
\caption{Illustration of operations $\lambda\pm\omega_l$ for
        $l=y_4=5$, applied to $\lambda=(8,6,5,3,3,1,1)$,
        translated to lattice paths.}
\label{fig:kleber-lattice2}
 \begin{picture}(25,8)(0,0)
  \linethickness{0.01pt}
  % \graphpaper[5](0,0)(25,19)
  \thinlines
  \put(4,0){
   \begin{picture}(15,6)(0,-1)
         % Heading
         \put(2,5){\hbox{Terminal points for
                        $(\lambda+\omega_5,\lambda-\omega_5)$:}}
         % Starting points, green
         \put(1,1){\circle{0.4}}
         \put(2,1){\circle{0.4}}
         \put(3,1){\circle{0.4}}
         \put(4,1){\circle{0.4}}
         \put(5,1){\circle{0.4}}
         \put(6,1){\circle{0.4}}
         \put(7,1){\circle{0.4}}
         % Starting points, blue
         \put(0,0.3){\circle*{0.4}}
         \put(1,0.3){\circle*{0.4}}
         \put(2,0.3){\circle*{0.4}}
         \put(3,0.3){\circle*{0.4}}
         \put(4,0.3){\circle*{0.4}}
         \put(5,0.3){\circle*{0.4}}
         \put(6,0.3){\circle*{0.4}}
         % End points, green
         \put(2,3.5){\circle{0.4}}
         \put(3,3.5){\circle{0.4}}
         \put(6,3.5){\circle{0.4}}
         \put(7,3.5){\circle*{0.4}}
         \put(10,3.5){\circle*{0.4}}
         \put(12,3.5){\circle*{0.4}}
         \put(15,3.5){\circle*{0.4}}
         % End points, blue
         \put(1,2.8){\circle*{0.4}}
         \put(2,2.8){\circle*{0.4}}
         \put(5,2.8){\circle{0.4}}
         \put(6,2.8){\circle*{0.4}}
         \put(9,2.8){\circle{0.4}}
         \put(11,2.8){\circle{0.4}}
         \put(14,2.8){\circle{0.4}}
         % Blocks
         \put(3.5,0.5){\oval(8,1.5)[b]} \put(6.5,0){\hbox{\bf 6}}
         \put(2,3.5){\oval(3,1.5)[t]}           \put(1.2,3.6){\hbox{\bf 5}}
         \put(6,3.5){\oval(3,1.5)[t]}           \put(5.2,3.6){\hbox{\bf 4}}
         \put(9.5,3.5){\oval(2,1.5)[t]} \put(9.2,3.6){\hbox{\bf 3}}
         \put(11.5,3.5){\oval(2,1.5)[t]}        \put(11.2,3.6){\hbox{\bf 2}}
         \put(14.5,3.5){\oval(2,1.5)[t]}        \put(14.2,3.6){\hbox{\bf 1}}
        \end{picture}
  }
 \end{picture}
\end{figure}

Likewise, adding some column of height $l=y_k$ to the ``blue object'' and
simultaneously removing such column from the ``green object'' amounts to interchanging colours of the
leftmost blue point and the rightmost green point in blocks $1,2,\dots,k$
if $k<n$; if $k=n$, then the same effect can be achieved by interchanging 
colours of the leftmost blue point and the rightmost green point in block
$n+1$. (See Figure~\ref{fig:kleber-lattice2}.)

\noindent
{\em Proof of Theorem~\ref{thm:Kleber}:\/}
Consider a two-coloured object from $s_\lambda s_\lambda$
in the lattice path interpretation.  As in Section~\ref{general}, we
look at the noncrossing perfect matching that the changing trails
induce among their $2n+2$ endpoints, the leftmost and rightmost points
in blocks $1,\ldots,n,n+1$.  Note that in the case of $s_\lambda
s_\lambda$, the parity constraint and the colour constraint of
Lemma~\ref{lem:parity} coincide.

Now consider the $k$ changing trails which begin at the leftmost
(blue) endpoints of blocks $1,2,\ldots,k$.  There are exactly two
% markus
cases:
% possibilities:

\begin{enumerate}
\item  
% markus 
% It is possible that the
The changing trails match these points up with the
rightmost (green) endpoints of blocks $1,2,\ldots,k$:
%  If this is the case, then
% markus
Then recolouring these $k$ trails results in an object of type
$s_{\lambda+\omega_l}s_{\lambda-\omega_l}$.  Conversely, given an
object of type $s_{\lambda+\omega_l}s_{\lambda-\omega_l}$, the parity
and colour constraints of Lemma~\ref{lem:parity} force the points
in blocks $1,2,\ldots,k$ to be matched amongst themselves, so they are
in bijection with this subset of $s_\lambda s_\lambda$ objects.

\item
Otherwise, some of those $k$ points must match up with points in
blocks $k+1,\ldots,n,n+1$.  Suppose there are $m$ such matchings, and
that they match the leftmost points in blocks $i_1<i_2<\cdots<i_m$
with points in blocks $j_1>j_2>\cdots>j_m$.  (Since the changing trails
cannot cross, we in fact know that $i_r$ is matched with $j_r$,
for $1\leq r\leq m$.)  Recolouring these $m$ trails gives an object of
type
$s_{\pi^{i_1\ldots i_m}_{j_1\ldots j_m}}(\lambda)
 s_{\mu^{i_1\ldots i_m}_{j_1\ldots j_m}}(\lambda)$.

This time, though, we do not have a bijection.  Given an object of 
$s_{\pi^{i_1\ldots i_m}_{j_1\ldots j_m}}(\lambda)
 s_{\mu^{i_1\ldots i_m}_{j_1\ldots j_m}}(\lambda)$,
the same parity and colour constraints of Lemma~\ref{lem:parity} do
guarantee that $m$ changing trails connect each $i_r$ with $j_r$.
However, when we recolour them to get an object of $s_\lambda
s_\lambda$, we may arrive at an object that has {\em other}
changing trails leaving blocks $1,\ldots,k$, aside from the $m$ we
% created.
% markus
considered.

% markus
Thus, we are in the ``typical'' situation for an inclusion--exclusion
argument, which immediately yields equation~\eqref{eq:kleber}.
% To resolve this dilemma, observe that if an object from $s_\lambda
% s_\lambda$ has $m$ changing trails leaving blocks $1,\ldots,k$ (as
% above), we could repeat the above construction recolour any subset of
% those trails.  Thus we use inclusion-exclusion to take the union of
% the preimages of the sets
% $s_{\pi^{i_1\ldots i_r}_{j_1\ldots j_r}}(\lambda)
% s_{\mu^{i_1\ldots i_r}_{j_1\ldots j_r}}(\lambda)$,
% and equation~\eqref{eq:kleber} follows.
\end{enumerate}

% markus
This finishes the proof.
\hfill\qed
\begin{rem}
When $k=1$, both cases of the above proof amount to recolouring the trail
beginning
at the rightmost green endpoint, so this is a special case of
Lemma~\ref{lem:general-skew}.  The $k=n$ case follows similarly, after
exchanging blue and green.
\end{rem}

%%%@@@ end input from parts/kleber

%%%@@@ input from parts/biblio
\ifx\undefined\bysame
\newcommand{\bysame}{\leavevmode\hbox to3em{\hrulefill}\,}
\fi

%%%@@@ end input from parts/biblio

\end{document}